\documentclass{elsarticle}
\usepackage{xcolor} 
\usepackage{amsfonts}
\usepackage{amssymb}
\usepackage{amsfonts}
\usepackage{mathtools}
\usepackage{enumitem}
\usepackage{color}

\usepackage{float}
\usepackage{tikz}
\usetikzlibrary{positioning,calc,arrows,automata}

\newtheorem{defn}{Definition}
\newtheorem{rmk}{Remark}
\newtheorem{pro}{Proposition}

\newcommand{\rd}{\textcolor{red}}

\newcommand{\bl}{\begin{list}{ \ }{
\leftmargin=.325in}}

\newcommand{\R}{\mathbb{R}} 
\newcommand{\CC}{\mathbb{C}} 
 
\newcommand{\VV}{\mathbb{V}} 
\newcommand{\WW}{\mathbb{W}}
\newcommand{\EE}{\mathbb{E}} 
\newcommand{\HH}{\mathbb{H}}

\newcommand{\A}{\mathcal{A}}
\newcommand{\B}{\mathcal{B}}
\newcommand{\C}{\mathcal{C}}
\newcommand{\E}{\mathcal{E}}
\newcommand{\I}{\mathcal{I}} 
\newcommand{\K}{\mathcal{K}}
\newcommand{\X}{\mathcal{X}}
\newcommand{\Y}{\mathcal{Y}}

\newcommand{\U}{\mathcal{U}}
\newcommand{\V}{\mathcal{V}}
\newcommand{\W}{\mathcal{W}}

\newcommand{\QQ}{\mathcal{Q}}

\newcommand{\fro}[1]{\left\vert\left\vert#1\right\vert\right\vert_F}

\usepackage{amssymb}
\usepackage{multirow}
\usepackage{ifthen}

\usepackage{float}
\usepackage[section]{placeins}
\usepackage{longtable}
\usepackage{mathrsfs}
\floatstyle{ruled}
\newfloat{algorithm}{htb}{alg}
\floatname{algorithm}{Algorithm}

\newcounter{algo@row}
\newcounter{algo@rowindent}
\newcommand{\algofont}[1]{\textbf{#1}}
\newcommand{\algonumbersize}[1]{\scriptsize{#1}}
\newcommand{\algopreitem}[1][\arabic{algo@row}]{\texttt{\algonumbersize{#1}}}
\newcommand{\algoitemskip}{\hspace{\value{algo@rowindent}cc}}

\newenvironment{algo}{\vskip.3em\small%
  \begin{list}{\algopreitem\texttt{\algonumbersize{:}}}{%
      \usecounter{algo@row}%
      \setcounter{algo@rowindent}{0}%
      \setlength{\itemindent}{2em}%
      \setlength{\labelwidth}{2em}
      \setlength{\parsep}{0cm}%
    }%
}{
  \end{list}\vskip-.5em
}
\newcommand{\algonewnestedopen}[2]{
  \newcommand{#1}[1][]{%
    \ifthenelse{\equal{##1}{}}{\item}{\item[{\algopreitem[##1]}]}
    \algoitemskip\algofont{#2}%
    \addtocounter{algo@rowindent}{1}%
    \ignorespaces
  }
}
\newcommand{\algonewnestedaux}[2]{
  \newcommand{#1}[1][]{
    \addtocounter{algo@rowindent}{-1}
    \ifthenelse{\equal{##1}{}}{\item}{\item[{\algopreitem[##1]}]}
    \algoitemskip\algofont{#2}%
    \addtocounter{algo@rowindent}{+1}%
    \ignorespaces
  }
}
\newcommand{\algonewnestedclose}[2]{
  \newcommand{#1}[1][]{
    \addtocounter{algo@rowindent}{-1}
    \ifthenelse{\equal{##1}{}}{\item}{\item[{\algopreitem[##1]}]}
    \algoitemskip\algofont{#2}%
    \ignorespaces
  }
}
\newcommand{\algonewcommand}[2]{
  \newcommand{#1}[1][default]{
    \ifthenelse{\equal{##1}{default}}{\item}{\item[{\algopreitem[##1]}]}%
    \algoitemskip\algofont{#2}%
    \ignorespaces
  }%
}
\newcommand{\algonewkeyword}[2]{\newcommand{#1}{\algofont{#2}}}
\algonewcommand{\STATE}{\ignorespaces}
\algonewcommand{\INPUT}{Input: }
\algonewcommand{\COMPUTE}{Compute: }
\algonewcommand{\BREAK}{Break }
\algonewcommand{\OUTPUT}{Output: }
\algonewnestedopen{\IF}{if }
\algonewnestedaux{\ELSEIF}{else if }
\algonewnestedaux{\ELSE}{else }
\algonewnestedclose{\ENDIF}{end if }
\algonewnestedopen{\FOR}{for }
\algonewnestedclose{\ENDFOR}{end for }
\algonewkeyword{\For}{for }%
\algonewkeyword{\To}{to }%
\algonewkeyword{\If}{if }%
\algonewkeyword{\Then}{then }%
\algonewkeyword{\Stop}{stop }%

\begin{document}           
\begin{frontmatter}
\title{A tensor formalism for multilayer network centrality measures using the
Einstein product}
\author{Smahane El-Halouy\fnref{addr1}}
\ead{elhalouysmahane@gmail.com}
\author{Silvia Noschese\fnref{addr2}}
\ead{noschese@mat.uniroma1.it}
\author{Lothar Reichel\fnref{addr3}}
\ead{reichel@math.kent.edu}
\address[addr1]{Department of Mathematical Sciences, Kent State University, Kent, 
OH 44242, USA, and Laboratory LAMAI, Faculty of Sciences and Technologies, Cadi Ayyad 
University, Marrakech, Morocco.}
\address[addr2]{Dipartimento di Matematica ``Guido Castelnuovo'', SAPIENZA Universit\`a 
di Roma, P.le A. Moro, 2, I-00185 Roma, Italy.}
\address[addr3]{Department of Mathematical Sciences, Kent State University, Kent, 
OH 44242, USA.}

\begin{abstract}
Complex systems that consist of diverse kinds of entities that interact in different 
ways can be modeled by multilayer networks. This paper uses the tensor formalism with the 
Einstein product to model this type of networks. Several centrality measures, that
are well known for single-layer networks, are extended to multilayer networks using 
tensors and their properties are investigated. In particular, subgraph centrality based 
on the exponential and resolvent of a tensor are considered. Krylov subspace methods based
on the tensor format are introduced for computing approximations of different measures for 
large multilayer networks. 
\end{abstract}

\end{frontmatter}

\section{Introduction}\label{sec1}
A network is a set of objects that are connected to each other in some fashion. 
Mathematically, a \emph{single-layer network} is represented by a graph $G=\{V,E\}$, where
the elements of the set $V=\{v_i\}_{i=1}^n$, referred to as vertices or nodes, represent 
the objects, and the elements of the set $E\subseteq V\times V$, designated as edges, 
represent the connections between the nodes. We denote an edge from node $v_i$ to node 
$v_j$ by $v_i\rightarrow v_j$.

Some real world examples require the modeling of more than one kind of nodes or of more 
than one type of edges. This holds, for instance, for the transportation network in a 
country when considering different means of transportation. The train and bus routes are 
different types of connections and should in some models be represented by different kinds
of edges. Moreover, train and bus stations may make up nodes with diverse properties. The 
connections between a train station and an adjacent bus station give rise to yet another 
kind of edges connecting different kinds of nodes, along which travelers typically walk. 
This kind of objects and connections 
can be modeled by \emph{multilayer networks}, which emphasize different kinds or 
connections, known as layers, between possibly different kinds of elements of a network. 
Each layer is represented by a single graph that contains the elements, or some of the 
elements, of the network and the connections between them in this layer. Edges connecting 
nodes from different layers model the interactions between different layers. Therefore, 
the nodes in a multilayer network require two indices, e.g., $v_i^\ell$, where the 
superscript $\ell$ denotes the layer, and the subscript $i$ determines the node in this 
layer. The set $V_L=V\times L$ represents all possible combinations of node-layers, where 
the set $V$ is made up of all nodes of the network considered. Each layer may be
made up of $V$ or some elements of $V$, and $L$ is the set of layers. The set of edges 
$E\subseteq V_L\times V_L$ represents all edges of the network. The 
special case when the set of nodes is the same in all layers, and edges that connect 
nodes in different layers are only allowed between a node and its copy in another layer,
is known as a \emph{multiplex network}. A nice recent paper by Bergermann and Stoll 
\cite{BS} studies multiplex networks and generalized matrix function-based centrality 
measures to this kind of networks. The authors use supra-adjacency matrices to represent 
multiplex networks. Recently, a global measure of communicability in a multiplex network, 
computed by means of the Perron root, and the right and left Perron vectors of the 
supra-adjacency matrix associated with this kind of network was introduced in 
\cite{ENR}. We are interested in using tensors for network analysis, because they
arise naturally when modeling multilayer networks.

\begin{figure}
\begin{center}
  \includegraphics[width=\textwidth]{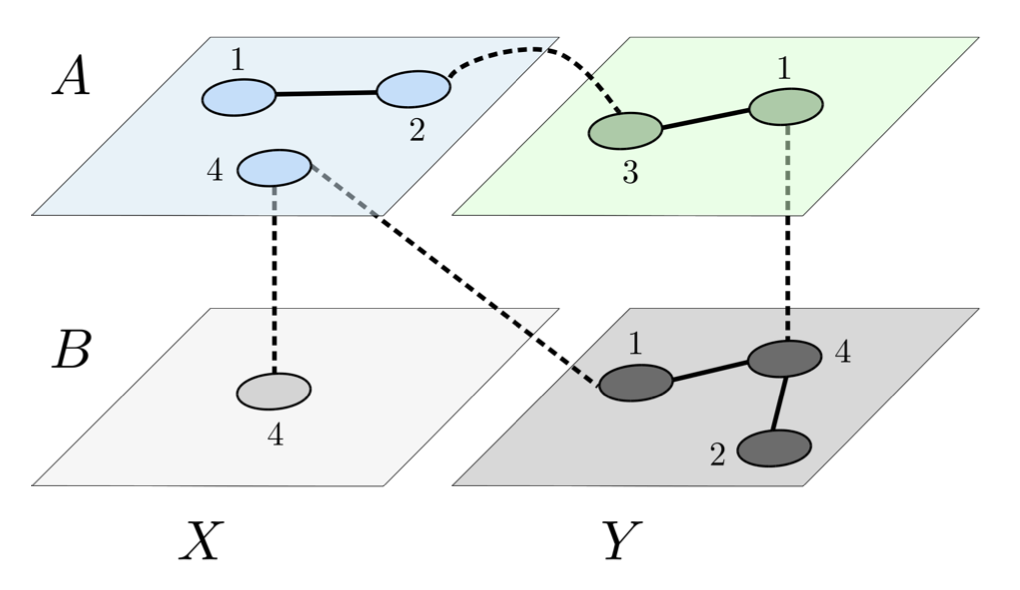}\
\end{center}
\caption{An example of a multilayer network with a set of four nodes $V=\{1,2,3,4\}$ and two aspects, the corresponding elementary layer sets $L_1=\{A,B\}$ and $L_2=\{X,Y\}$. The total number of layers is four and they are $(A,X)$, $(A,Y)$, $(B,X)$ and $(B,Y)$. Each layer includes some of the elements of $V$.}\label{figEx}
\end{figure}

The model mentioned above can be generalized to represent not only networks with multiple 
layers but also different aspects. To allow for the modeling of more than one aspect, we 
define a sequence $\{L_j\}_{j=1}^{d}$ of sets of elementary layers with $d$ being the 
number of aspects that we would like to model; $L_j$ is the set of layers for aspect $j$.
Then the total number of layers is $|L_1|\times |L_2|\times\ldots\times |L_d|$ and we have
$V_L=V\times L_1\times \ldots \times L_d$. The nodes now are identified by using $d+1
$ indices $v_i^{\ell_1,\ldots ,\ell_d}$, where the subscript $i$ indicates the number of 
the node and the superscript $\ell_1,\ldots,\ell_d$ shows the specific layer. For more 
details on this kind of generalization, we refer to \cite{DSC,KBGMP} and the references 
therein, where general frameworks for multilayer network are discussed together with their 
mathematical formulation. Figure \ref{figEx} illustrates a simple multilayer network with 
$2$ aspects; this figure can also be found in
\cite{BP}. An example of a real multilayer network with multiple aspects in biology is
provided in \cite{MGMOP}, where the first aspect is the type of data (genomic, metabolomic,
or proteomic), and the second aspect models different biological pathways; see Figure 2 in 
\cite{MGMOP}.

Single-layer networks are often represented by an adjacency matrix, which is helpful for
extracting information about the network, e.g, by evaluating functions of the the adjacency
matrix or by computing certain eigenvectors of this matrix. For instance, Estrada and Higham 
\cite{EH2} describe how the matrix 
exponential and resolvent can be used to determine how easy it is to communicate between 
nodes in a single-layer network, and which nodes are the most important ones; see also 
Estrada \cite{Esbook} and references therein. For multilayer networks, we use tensors, 
i.e., a multidimensional generalization of matrices, and represent the network by a 
$2(d+1)$-order adjacency tensor $\A$ of size 
$(|V|\times|L_1|\times\ldots\times|L_d|)\times(|V|\times|L_1|\times\ldots\times|L_d|)$,
where $|V|$ denotes the total number of nodes, and $|L_j|$ 
designates the number of layers for property $j$, for $j=1,2,\ldots,d$. The entry 
$\A(i,\ell_1,\ldots,\ell_d,j,k_1,\ldots,k_d)$ of the tensor $\A$ for an \emph{unweighted}
multilayer network with sets of layers $L_j$, $j=1,2,\ldots,d$, is one if there is an
edge $v_i^{\ell_1,\ldots,\ell_d}\rightarrow v_j^{k_1,\ldots,k_d}$; otherwise the tensor 
entry is zero. For a \emph{weighted} network a tensor entry $1$ may be replaced by a real, 
generally positive, number. In a directed network some of the edges represent ``one-way
streets''. Note that some nodes may not be present in all layers. Therefore, considering 
empty nodes is necessary to allow the tensorial representation. For instance, the network 
illustrated in Figure \ref{figEx} can be represented by a $6^{th}$ order tensor of size 
$4\times 2\times 2\times 4\times 2\times 2$ by adding empty nodes so that every layer is 
made up of $4$ nodes.

Adjacency tensors allow us to capture the structure and complexity of relationships 
between the nodes of a multilayer network. We are interested in investigating and 
generalizing some centrality measures that are well established for single-layer
networks to multilayer networks by using the tensor formalism and applying tensor tools,
such as the Einstein product and tensor functions. Several centrality measures have been 
studied for multilayer networks using the tensor formalism in \cite{DSOGA}. Eigenvector 
multicentrality has been investigated for multilayer networks via a tensor-based framework
in \cite{WZCSLZP}; see also \cite{CRZT}. In addition to generalizing centrality measures 
that are commonly used for single-layer networks to multilayer networks, we describe 
practical and efficient ways to compute these measures by using Krylov subspace methods 
based on the tensor format.

This paper is organized as follows. Tensor notation, definitions, and properties used 
throughout this paper are described in Section \ref{sec2}. Section \ref{sec3} discusses 
the extension of matrix functions to tensor functions using the Einstein tensor 
product. We define centrality measures for multilayer networks based on the tensor 
representation and tensor functions. Section \ref{sec4} describes Krylov subspace methods 
based on the tensor format using the Einstein product, and discusses their application to 
the approximation of tensor functions. Section \ref{sec5} presents a few computed examples
and Section \ref{sec6} contains concluding remarks.

\section{Preliminaries}\label{sec2}
This section presents notation and properties of tensors that will be used throughout this
paper. We start with a generalization of the matrix-matrix product to tensors that is 
referred to as the Einstein product.

\begin{defn}[Einstein Product]
Let $\A\in\R^{I_1\times\ldots \times I_N\times  J_1\times\ldots \times J_M}$ and 
$\B\in\R^{J_1\times\ldots \times J_M\times K_1\times\ldots \times K_L}$ be tensors of 
orders $N+M$ and $M+L$, respectively. The product 
$\C=\A\ast_M\B\in\R^{I_1\times\ldots \times I_N\times K_1\times\ldots \times K_L}$ of the
tensors $\A$ and $\B$ is a tensor of order $N+L$ with entries
\[
\C_{i_1,\ldots ,i_N,k_1,\ldots ,k_L}=\sum_{j_1,\ldots,j_M}\A_{i_1,\ldots,i_N,j_1,\ldots,
j_M}\B_{j_1,\ldots ,j_M,k_1,\ldots ,k_L}.
\]
It is commonly referred to as the Einstein product; see \cite{BM,BLNT,GIJB}. The subscript 
$M$ in $\ast_M$ indicates the last and first $M$ dimensions of $\A$ and $\B$, respectively, 
over which the sum is evaluated. 
\end{defn}

The identity tensor $\I=[\I_{i_1,\ldots,i_N,j_1,\ldots,j_N}]\in
\R^{I_1\times\ldots\times I_N\times I_1\times\ldots\times I_N}$ under the Einstein product 
has the entries
\[
\I_{i_1,\ldots,i_N,j_1,\ldots,j_N}=\left\{\begin{array}{cl} 1, &
~{\rm if}~i_k=j_k, ~{\rm for}~k=1,2,\ldots,N, \\ 0, &~{\rm otherwise}.\end{array}\right.
\]
\begin{rmk}
A tensor of even order 
$\A\in\R^{I_1\times\cdots\times I_N\times J_1\times\cdots\times J_N}$ is said to be 
\emph{square} if the first set of dimensions equals the second set, i.e., if $I_k=J_k$ for
$k=1,2,\ldots,N$; see, e.g., \cite{QL}. The adjacency tensor of a multilayer network is of
even order and square, however, for some computations tensors of different orders are 
required. This is possible when using the Einstein product by choosing a suitable number of 
dimensions over which we carry out the summation.
\end{rmk}

\begin{rmk}
The transpose of a tensor 
$\A\in\R^{I_1\times\ldots \times I_N\times  J_1\times\ldots \times J_M}$ is a tensor 
$\B\in\R^{J_1\times\ldots \times J_M\times I_1\times\ldots \times I_N}$  such that 
$\A_{i_1,\ldots ,i_N,j_1,\ldots ,j_M}=\B_{j_1,\ldots ,j_M,i_1,\ldots ,i_N}$; see, e.g.,
\cite{QL}. 

The $n$-mode product is a well-known tensor-matrix product; see \cite{KB2}. For an 
$N^{th}$ order tensor $\A\in\R^{I_1\times\ldots \times I_N}$ and a matrix 
$A\in\R^{I_n\times J}$, their $n$-mode product is an $N^{th}$ order tensor 
$\A\times_nA\in\R^{I_1\times\ldots I_{n-1}\times J\times I_{n+1} \times\ldots I_N}$.

If $\A\in\R^{I_1\times\ldots \times I_N\times J}$ is an $(N+1)^{th}$ order tensor and $A\in\R^{J\times I}$, then the $n$-mode product of $\A$ and $A$ over mode $N+1$ is the same as the Einstein product when summing over the last mode of 
the tensor. In other words, we have 
\[
\A\ast_1A = \A\times_{N+1}A^T.
\]
\end{rmk}

We can reorganize the entries of a tensor in different ways to obtain a $2D$ array, i.e.,
a matrix. This transformation is known as \emph{matricization} or \emph{flattening}. 
We flatten a tensor by using lexicographical ordering of the indices.

\begin{defn}[Tensor flattening]\label{Def2}
Let $\A\in\R^{I_1\times\ldots \times I_N\times  J_1\times\ldots \times J_M}$ be a tensor
of order $N+M$. The elements of the matrix $A\in\R^{I_1\ldots I_N\times J_1\ldots J_M}$
obtained by flattening the tensor $\A$ are given by
\[
A_{i,j}=\A_{i_1,\ldots ,i_N,j_1,\ldots ,j_M},
\]
where 
\begin{eqnarray*}
i&=&i_1+\sum_{p=2}^N(i_p-1)\prod_{q=1}^{p-1}I_q,\\
j&=&j_1+\sum_{p=2}^M(j_p-1)\prod_{q=1}^{p-1}J_q.
\end{eqnarray*}
Here and below the indices $i_p$ and $j_p$ live in their domains, i.e., 
$1\leq i_p\leq I_p$ for $1\leq p\leq N$, and $1\leq j_q\leq J_p$ for $1\leq q\leq M$.
We define $A={\rm mat}(\A)$ and $\A={\rm mat}^{-1}(A)$.
\end{defn}

\begin{rmk}
For a multiplex network, ${\rm mat}(\A)$ is the supra-adjacency matrix defined in 
\cite{BS}.
\end{rmk}

\begin{pro}\label{prop1}
Let $\A\in\R^{I_1\times\ldots\times I_N\times J_1\times\ldots\times J_M}$ and 
$\B\in\R^{J_1\times\ldots\times J_M\times K_1\times\ldots\times K_L}$ be tensors of orders
$N+M$ and $M+L$, respectively. Then
\[
{\rm mat}(\A\ast_M\B)={\rm mat}(\A)\cdot{\rm mat}(\B),
\]
where $\cdot$ denotes the usual matrix product.
\end{pro}

\noindent
Proof: 
The result follows by direct computations; see \cite{BLNT} for details.~~~$\Box$

\begin{defn}\label{def3}
The following definitions can be found in, e.g., \cite{BLNT,KB1,QL}.
\begin{enumerate}
\item The trace of a square tensor, 
$\A\in\R^{I_1\times\ldots\times I_N\times I_1\times\ldots\times I_N}$, is given by 
\[
{\rm tr}(\A)=\sum_{i_1,\ldots ,i_N}\A_{i_1,\ldots ,i_N,i_1,\ldots ,i_N}.
\]
\item The inner product of two tensors of the same size
$\X,\Y\in\R^{I_1\times\ldots\times I_N\times J_1\times\ldots\times J_M}$ is defined as 
\[
\langle\X,\Y\rangle=\sum_{i_1,\ldots,i_N,j_1,\ldots,j_M}\X_{i_1,\ldots,i_N,j_1,\ldots,j_M}
\Y_{i_1,\ldots,i_N,j_1,\ldots,j_M}.
\]
For square tensors, we have 
\[
\langle\X,\Y\rangle={\rm tr}(\X^T\ast_N\Y).
\]
\item The Frobenius norm of a tensor is given by 
\[
\fro{\X}=\sqrt{ \langle\X,\X\rangle}.
\]
If $\X$ is a square tensor, then
\[
\fro{\X}=\sqrt{{\rm tr}(\X^T\ast_N\X)}.
\]
\item For a positive integer $p$, we define the $p^{\rm th}$ power of a square tensor 
$\A\in\R^{I_1\times\ldots\times I_N\times I_1\times\ldots\times I_N}$ by using the Einstein 
product recursively as
\[
\A^p=\A\ast_N\A^{p-1},
\]
where $\A^0=\I$ is the identity tensor.
\end{enumerate}
\end{defn}

\begin{pro}\label{prop2}
One has
\begin{equation}\label{prop2(1)}
\fro{\A}=\fro{{\rm mat}(\A)}
\end{equation}
and, if $\A$ is a square tensor, then
\begin{equation}\label{prop2(2)}
\fro{\A^p}\leq\fro{\A}^p.
\end{equation}
\end{pro}

\noindent
Proof: Let $\A\in\R^{I_1\times\ldots\times I_N\times J_1\times\ldots\times J_M}$ be a 
tensor of order $N+M$. One has
\[
\fro{\A}=\sqrt{ \langle\A,\A\rangle}=\sqrt{\sum_{i_1,\ldots,i_N,j_1,\ldots,j_M}\A_{i_1,
\ldots,i_N,j_1,\ldots,j_M}^2}=\sqrt{\sum_{i=1}^N\sum_{j=1}^M A_{i,j}^2}=\fro{A},
\]
with $A={\rm mat}(\A)$ given in Definition \ref{Def2}. This shows \eqref{prop2(1)}. Assume
now that $M=N$. According to Definition \ref{def3}({\it 4}.), 
by applying Proposition \ref{prop1} $p-1$ times, one has ${\rm mat}(\A^p)=({\rm mat}(\A))^p$,
so that $\fro{{\rm mat}(\A^p)}\leq \fro{{\rm mat}(\A)}^p$, where the inequality is due to the 
submultiplicativity of the Frobenius matrix norm. Thus, thanks to \eqref{prop2(1)}, 
one has inequality \eqref{prop2(2)}. ~~~$\Box$ 
\begin{rmk}\label{rmk3}
The eigenvalue problem for $4^{\rm th}$-order tensors is discussed in \cite{DSC}. We can express this problem for any square tensor by using the Einstein product; we have
\[
\A\ast_N \X = \lambda \X,
\]
where $\A\in\R^{I_1\times\ldots \times I_N\times I_1\times\ldots\times I_N}$, 
$X\in\R^{I_1\times\ldots\times I_N}$ and $\lambda\in\CC$. This eigenvalue problem is 
equivalent to the matrix eigenvalue problem 
${\rm mat}(\A)\cdot {\rm vec}(\X) = \lambda {\rm vec}(\X)$, where {\rm vec} is the 
vectorization operator that takes a tensor and rearranges it into a single column vector by 
concatenating its elements, i.e., it stacks the elements $\X$ to form a column vector.
\end{rmk}

\section{Tensor functions and centrality measures for multilayer networks}\label{sec3}
Centrality measures have been thoroughly studied for single-layer networks. These measures
include Katz centrality, subgraph centrality, and total communicability with respect to a 
node; see, e.g., \cite{BK,DMR,Esbook,EH1,EH2,ERV}. In this section, we introduce analogues 
of these measures for multilayer networks using the adjacency tensor and the Einstein 
product. 

We recall that, for single-layer networks, a walk from node $v_{i_1}$ to node $v_{i_j}$ is
defined as a sequence of edges 
\[
v_{i_1}\rightarrow v_{i_2},~v_{i_2}\rightarrow v_{i_3},~\ldots,~
v_{i_{j-1}}\rightarrow v_{i_j}
\]
that can be traversed to reach node $v_{i_j}$ from node $v_{i_1}$. The length of the walk 
is the number of edges, $j-1$. Nodes and edges may be repeated in a walk. A walk is said
to be short if $j-1$ is fairly small. There may be an edge from node $i$ in layer 
$(\ell_1,\ell_2,\ldots ,\ell_d)$ to node $j$ in layer $(k_1,k_2,\ldots ,k_d)$ in a 
multilayer network. We denote this edge by 
$v_i^{\ell_1,\ell_2,\ldots ,\ell_d}\rightarrow v_j^{k_1,k_2,\ldots ,k_d}$. A walk in a multilayer 
network from node $v_{i_1}^{\ell_{k_1},\ldots ,\ell_{k_d}}$ to node $v_{i_j}^{\ell_{k_1+j},\ldots ,\ell_{k_d
+j}}$ is defined as a sequence of edges such as
\begin{eqnarray*}
&&v_{i_1}^{\ell_{k_1},\ldots,\ell_{k_d}}\rightarrow 
v_{i_2}^{\ell_{k_1+1},\ldots,\ell_{k_d+1}},~~v_{i_2}^{\ell_{k_1+1},\ldots ,\ell_{k_d+1}}
\rightarrow v_{i_3}^{\ell_{k_1+2},\ldots,\ell_{k_d+2}},\ldots,\\
&&v_{i_{_j-1}}^{\ell_{k_1+j-1},\ldots ,\ell_{k_d+j-1}}\rightarrow 
v_{i_j}^{\ell_{k_1+j},\ldots ,\ell_{k_d+j}}
\end{eqnarray*}
that can be traversed to reach node $v_{i_j}^{\ell_{k_1+j},\ldots ,\ell_{k_d+j}}$ from node $v_{i_1}^{\ell_{k_1},\ldots ,\ell_{k_d}}$. The length
of the walk is the number of edges, $j-1$, and nodes and edges may be repeated in a walk. 
A closed multilayer walk is a multilayer walk for which the starting and ending nodes are 
the same, i.e., $i_1=i_j$ and $\ell_{k_n}=\ell_{k_n+j}$, for $n=1,2,\ldots ,d$. Estrada \cite{E} has defined a walk in a 
multiplex network in an analogous fashion. 
 
The entries of the adjacency tensor 
$\A\in\R^{N\times K_1\times\ldots K_d\times N\times K_1\times\ldots K_d}$ of an
unweighted undirected multilayer
network tell us whether there is an edge between any pair of nodes (between the 
same or different layers).  The entries of the Einstein product of the adjacency tensor 
$\A$ with itself, 
\begin{eqnarray*}
\B_{i,\ell_1,\ell_2,\ldots ,\ell_d,j,k_1,k_2,\ldots ,k_d}&=&\A\ast_{d+1}\A\\
&=&\sum_{\substack{p,q_1,\ldots ,q_d}}\A_{i,\ell_1,\ell_2,\ldots ,\ell_d,p,q_1,\ldots ,q_d}\A_{p,q_1,\ldots ,q_d,j,k_1,k_2,\ldots ,k_d},
\end{eqnarray*}
where $1\leq i,j\leq N,~~1\leq \ell_s,k_s\leq L,~~1\leq s\leq d$, show the number of multilayer walks of length $2$ between pairs of nodes $v_i^{\ell_1,\ell_2,\ldots ,\ell_d}$ and 
$v_j^{k_1,k_2,\ldots ,k_d}$. Similarly, let $p$ be a positive integer. Then the entries of the tensor
$\A^p=\A\ast_{d+1}\A^{p-1}$ display the number of multilayer walks of length $p$ between 
pairs of nodes $v_i^{\ell_1,\ell_2,\ldots ,\ell_d}$ and $v_j^{k_1,k_2,\ldots ,k_d}$. In addition, the entry $\A^p_{i,\ell_1,\ell_2,\ldots ,\ell_d,i,\ell_1,\ell_2,\ldots ,\ell_d}$ 
provides the number of closed multilayer walks of length $p$ that start at node 
$v_i^{\ell_1,\ell_2,\ldots ,\ell_d}$. This suggests, in order to take into account walks 
of all possible lengths $p\geq 0$, the introduction of the tensor function
\begin{equation}\label{tenfun}
f(\A) = \sum_{p=0}^{\infty}c_p\A^p,
\end{equation}
where the coefficients $c_p$ generally are real and nonnegative, and are chosen so that
the series converges.

Let the tensor $\E_{s,t_1,\ldots ,t_d}\in\R^{N\times K_1\times\ldots \times K_d}$ with all
entries equal to zero except for the $(s,t_1,\ldots ,t_d)^{\rm th}$ entry, which is one. 
We refer to the entry 
\begin{equation}\label{com}
f(\A)_{i,\ell_1,\ell_2,\ldots ,\ell_d,j,k_1,k_2,\ldots,k_d}=
\E_{i,\ell_1,\ell_2,\ldots,\ell_d}\ast_{d+1} f(\A) \ast_{d+1} \E_{j,k_1,k_2,\ldots ,k_d}
\end{equation}
as the \emph{communicability} from node $v_i^{\ell_1,\ell_2,\ldots ,\ell_d}$ to node 
$v_j^{k_1,k_2,\ldots ,k_d}$; a relatively large 
value indicates that it is easy to send information from node $v_i^{\ell_1,\ell_2,\ldots ,\ell_d}$ to node $v_j^{k_1,k_2,\ldots ,k_d}$. 
Moreover, we refer to the entry 
\begin{equation}\label{sgc}
f(\A)_{i,{\ell_1,\ell_2,\ldots ,\ell_d},i,{\ell_1,\ell_2,\ldots ,\ell_d}}=\E_{i,{\ell_1,\ell_2,\ldots ,\ell_d}}\ast_{d+1} f(\A) \ast_{d+1} \E_{i,{\ell_1,\ell_2,\ldots ,\ell_d}}
\end{equation}
as the \emph{subgraph centrality} of node $v_i^{\ell_1,\ell_2,\ldots ,\ell_d}$; a relatively large value indicates that
this node is important because much information may pass through it. These notions of 
communicability and subgraph centrality for nodes in multilayer graphs generalize the 
analogous definitions introduced and explored by Estrada and his collaborators, as well as
others, in \cite{ADR,DMR,EH1,EH2,ERV} for single-layer networks. In these references $\A$ 
in \eqref{tenfun} is replaced by the adjacency matrix for the single-layer graph and, 
hence, $f$ is a matrix function. 

In many network applications short walks are more important than long walks, because it is 
easier to transmit information via a few edges than via many edges. This suggests that the
coefficients $c_p$ should satisfy $0\leq c_{p+1}\leq c_p$ for large $p$-values. One of the 
most commonly used matrix functions for single-layer networks is the matrix exponential. 
For multilayer networks, we therefore introduce the tensor exponential
\[
\exp(\beta\A) = \sum_{p=0}^{\infty}\frac{\beta^p\A^p}{p!}.
\]
It follows from Propositions \ref{prop1} and \ref{prop2} that this series converges for 
any fixed $\beta$ in the interval $0\leq\beta<\infty$. Since the first term $\I$ has no 
natural interpretation in the context of network modeling, we will use the modified tensor
exponential 
\begin{equation}\label{modexp}
\exp_0(\beta\A) := \exp(\beta\A) -\I,
\end{equation}
where $\I$ denotes the identity tensor.

We introduce the $\delta$-\emph{effective diameter} of the network determined by the 
tensor function \eqref{tenfun}. It is defined as the smallest integer $k$ such that
\begin{equation}\label{effdiam}
\frac{\max_{\ell>k}|c_{\ell}|}{\max_{1\leq j\leq k}|c_j|}\leq\delta,
\end{equation}
for some small $\delta>0$. Roughly, this diameter is $k$ if the tensor function 
\eqref{tenfun} can be approximated well by a polynomial of degree $k$. This means that 
walks of length larger than $k$ do not significantly affect the properties of the network. 
for some small $\delta>0$. Its importance for the communicability in single-layer networks
is explored in \cite{ADR}. The definition of the effective diameter in \cite{ADR} differs 
slightly from \eqref{effdiam} and is for matrix functions.

\begin{pro}
Let the tensor function \eqref{tenfun} be the modified tensor exponential \eqref{modexp}
for some $\beta>0$. Then the left-hand side of \eqref{effdiam} decreases to zero as $k$
increases.
\end{pro}

\noindent
Proof: Assume that $\beta>1$ and let $k_\beta$ denote the integer part of $\beta$. Then
\[
\max_{\ell>k} c_\ell=\left\{\begin{array}{cc}
\displaystyle{\frac{\beta^{k_\beta}}{k_\beta!}} & \mbox{~~if~~} k<k_\beta,\\[3mm]
\displaystyle{\frac{\beta^{k+1}}{(k+1)!}} & \mbox{~~if~~} k\geq k_\beta,
\end{array}\right.
\]
and
\[
\max_{1\leq j\leq k} c_j=\left\{\begin{array}{cc}
\displaystyle{\frac{\beta^{k_\beta}}{k_\beta!}} & \mbox{~~if~~} k\geq k_\beta,\\[3mm]
\displaystyle{\frac{\beta^k}{k!}} & \mbox{~~if~~} k < k_\beta.
\end{array}\right.
\]
It follows that
\[
\frac{\max_{\ell>k} c_\ell}{\max_{1\leq j\leq k}c_j}=\left\{\begin{array}{cc}
\displaystyle{\frac{\beta^{k_\beta}}{k_\beta!}\cdot\frac{k!}{\beta^k}} & \mbox{~~if~~} k<k_\beta, \\[3mm]
\displaystyle{\frac{\beta^{k+1}}{(k+1)!}\cdot\frac{k_\beta!}{\beta^{k_\beta}}} & \mbox{~~if~~} k\geq k_\beta.
\end{array}\right.
\]
Therefore this quotient converges to zero as $k$ increases.

We turn to the situation when $0<\beta\leq 1$. Then $\max_{1\leq j\leq k}c_j=c_1$ and 
\[
\max_{\ell>k}c_\ell=c_{k+1}=\frac{\beta^{k+1}}{(k+1)!}\to 0 \mbox{~~as~~} k\to\infty,
\]
and the proposition follows. $~~~\Box$

Resolvents of the adjacency matrix also are commonly used to determine properties of nodes
in a single-layer network; see, e.g., Estrada and Higham \cite{EH2} and Katz \cite{Ka}. We 
define the modified tensor resolvent,
\begin{equation}\label{modres}
{\rm res}_0(\A,\alpha)=(\I-\alpha\A)^{-1}-\I = \sum_{p=1}^{\infty}\alpha^p\A^p,
\end{equation}
which is convergent for $0<\alpha<1/|\lambda_{\max}|$, where $\lambda_{\max}$ denotes an
eigenvalue of largest magnitude of $\A$. For many adjacency tensors of interest, 
$\lambda_{\max}$ is real and positive. Conditions under which this is the case are 
discussed by Qi and Luo \cite{QL}. The eigenvalue 
$\lambda_{\max}$ can be computed as an eigenvalue of a matrix using the relation of Remark
\ref{rmk3}. The choice of $\alpha$ affects the $\delta$-effective diameter of the tensor 
function \eqref{modres}. This is discussed for the matrix resolvent in \cite{ADR}.

We also define the \emph{multilayer total communicability of node} $v_i^{\ell_1,\ell_2,\ldots ,\ell_d}$ by 
\begin{equation}\label{mltciell}
\E_{i,\ell_1,\ell_2,\ldots ,\ell_d}\ast_{d+1} f(\A) \ast_{d+1} \E,
\end{equation}
where  $\E\in\R^{N\times K_1\times\ldots\times K_d}$ is a tensor with all entries equal 
to one, and the \emph{multilayer total communicability} by
\begin{equation}\label{mltc}
\E\ast_{d+1} f(\A) \ast_{d+1} \E.
\end{equation}
The latter definitions generalize analogous notions introduced for single-layer networks
by Benzi and Klymko \cite{BK} and Katz \cite{Ka} for $f$ being the matrix exponential or a
matrix resolvent. We refer to the quantity defined in (\ref{mltciell}) as the \emph{multilayer Katz centrality} when $f$ is the tensor resolvent. 

For small networks, we can evaluate the tensor functions discussed above by applying the 
flattening operator ${\rm mat}$, its inverse ${\rm mat}^{-1}$, Proposition \ref{prop1}, 
and using the following result.

\begin{pro}\label{prop3}
Let the tensor function $f$ be defined by \eqref{tenfun} with a power series that
converges sufficiently rapidly. Then
\[
f(\A) = {\rm mat}^{-1}(f({\rm mat}(\A))).
\]
\end{pro}

\noindent
Proof: By Proposition \ref{prop1}, ${\rm mat}(\A^p)=({\rm mat}(\A))^p$. Hence, 
$\sum_{p=1}^n{\rm mat}(c_p\A^p) =\sum_{p=1}^nc_p({\rm mat}(\A))^p$ and we have

\[
f({\rm mat}(\A))=\lim_{n\to\infty} \sum_{p=1}^nc_p({\rm mat}(\A))^p=
\lim_{n\to\infty}{\rm mat}(  \sum_{p=1}^nc_p\A^p), 
\]
and by the definition of ${\rm mat}$, we can write 
\[
\lim_{n\to\infty}{\rm mat}( \sum_{p=1}^nc_p\A^p)={\rm mat}(\lim_{n\to\infty} 
\sum_{p=1}^nc_p\A^p).
\]
 
Thus, one has
\[
{\rm mat}(f(\A))=f({\rm mat}(\A)).
\]
Applying the inverse operator ${\rm mat}^{-1}$ to both sides concludes the proof.
$~~~\Box$

The evaluation of tensor functions using the above proposition is feasible for tensors 
that represent small to medium-sized multilayer networks. However, the computations are 
very demanding for large-scale multilayer networks. Approximations of tensor functions
for the latter kind of networks can be computed fairly inexpensively by applying Krylov
subspace methods to the flattened adjacency tensor, i.e., supra-adjacency matrix, as in 
\cite{BS}. However, our main goal is to contribute to the development of a formalism where
tensors are used. Therefore, we suggest computing communicability and centrality measures 
using Krylov subspace methods based on the tensor format. This is discussed in the 
following section.

\section{Krylov subspace methods}\label{sec4}
Krylov subspace methods are well suited to approximate many matrix functions; see, e.g.,
\cite{BR,GM} for illustrations. They also have been applied successfully to the 
approximation of tensor functions and the solution of tensor systems of equations; see
\cite{BJNR,BGJR,GIJB,HYM,Lu} and references therein. It is therefore natural to seek to 
approximate the tensor functions mentioned in the previous section by Krylov subspace 
methods. We first discuss the application of the global tensor Arnoldi process to the 
approximation of multilayer centrality measures, and subsequently consider the tensor
block Arnoldi process.

\subsection{A global tensor Arnoldi process based on the Einstein product}
The global matrix Arnoldi process is a Krylov subspace method that was introduced by 
Jbilou et al. \cite{JMS,JST} for the reduction of a large matrix to a small one. A global
tensor Arnoldi process for the reduction of a large tensor to a small matrix using the 
Einstein product is described by El Guide et al. \cite{GIJB}. The application of $m$ steps
of the latter process to the tensor $\A\in\R^{N\times K_1\times\ldots\times K_d\times N\times K_1\times\ldots\times K_d }$ with initial tensor 
$\V\in\R^{N\times K_1\times\ldots\times K_d}$ determines, when no breakdown occurs, an 
orthonormal basis for 
the tensor Krylov subspace
\begin{equation}\label{KrSub}
\K_{m+1}(\A,\V)={\rm span}\{\V,\A\ast_{d+1} \V,\ldots ,\A^m\ast_{d+1} \V\}:=
\left\{\sum_{i=0}^m\omega_i\A^i\ast_{d+1} \V,~ \omega_i\in\R\right\}.
\end{equation}
This definition of the subspace is analogous to the definition of the
solution subspace for global matrix methods used in \cite{JST}.

\begin{algorithm}
\caption{Global tensor Arnoldi process}\label{Algo1}
\begin{algo}
\INPUT Adjacency tensor $\A\in\R^{N\times K_1\times\ldots\times K_d\times N\times K_1\times\ldots\times K_d }$, initial tensor 
$\V\in\R^{N\times K_1\times\ldots\times K_d}$, and number of steps $m$.
\OUTPUT Orthonormal basis $\VV_{m+1}=\{\V_1,\V_2,\ldots,\V_{m+1}\}$ for the tensor Krylov
subspace \eqref{KrSub} and nontrivial entries of the upper Hessenberg matrix 
$H_{m+1,m}=[h_{ij}]\in\R^{(m+1)\times m}$.
\STATE $\V_1=\V/\|\V\|_F$
\FOR $j=1,\ldots,m$ do
\STATE $\W = \A\ast_{d+1}\V_j$ 
\FOR $i=1,\ldots ,j$ do 
\STATE $h_{i,j}=\langle \V_i,\W\rangle$ 
\STATE $\W=\W-h_{ij}\V_j$
\ENDFOR
\STATE $h_{j+1,j}=\fro{\W}$ 
\IF $h_{j+1,j}=0$, then \Stop
\ELSE $\V_{j+1}=\W/h_{j+1,j}$
\ENDIF
\ENDFOR
\end{algo}
\end{algorithm} 

The computations are described by Algorithm \ref{Algo1}. The algorithm is said to 
\emph{break down} at step $j$ if $h_{i+1,i}>0$ for $1\leq i<j$ and $h_{j+1,j}=0$. In the
absence of breakdown, the algorithm determines the tensor
$\VV_{m+1}=[\V_1,\V_2,\ldots,\V_{m+1}]\in\R^{N\times K_1\times\ldots\times K_d\times(m+1)}$ with orthonormal block 
columns, i.e.,
\[
\langle \V_i,\V_j\rangle:=\V_i\ast_{d+1} \V_j=
\left\{\begin{array}{cc} 1, & i= j,\\ 0, & i \ne j,\end{array}\right.
\]
that span the tensor Krylov subspace $\K_{m+1}(\A,\V)$. In line 5 of Algorithm \ref{Algo1},
we have
\[
\A\ast_{d+1}\VV_j=[\A\ast_{d+1}\V_1,\A\ast_{d+1}\V_2,\ldots,\A\ast_{d+1}\V_j]\in\R^{N\times K_1\times\ldots\times K_d\times j}.
\]
It follows from the recursion relation of Algorithm \ref{Algo1} that
\begin{equation}\label{Glob1}
\A\ast_{d+1}\VV_m=\VV_{m+1}\ast_1 H_{m+1,m},
\end{equation}
where $H_{m+1,m}=[h_{ij}]\in\R^{(m+1)\times m}$ is an upper Hessenberg matrix made up of 
the coefficients $h_{ij}$ generated in lines 7 and 10 of Algorithm \ref{Algo1}; all entries 
below the subdiagonal of $H_{m+1,m}$ vanish. 

Let the matrix $H_m\in\R^{m\times m}$ be obtained by deleting the last row of $H_{m+1,m}$.
Then 
\[
\VV_m^T\ast_{d+1}\A\ast_{d+1}\VV_m = H_m,
\]
where $\VV_m=[\V_1,\V_2,\ldots,\V_m]\in\R^{N\times K_1\times\ldots\times K_d\times m}$ and 
$\VV_m^T=[\V_1,\V_2,\ldots,\V_m]^T\in\R^{m\times N\times K_1\times\ldots\times K_d}$. Hence, $H_m$ is the orthogonal 
projection of $\A$ onto the subspace $\K_m(\A,\V)$ with respect the basis $\VV_m$. This 
suggest to use the approximation
\begin{equation}\label{fred}
\VV_m\ast_1f(H_m)\ast_1E_1 \|\V\|_F 
\end{equation}
of $f(\A)\ast_{d+1} \V$, where $\|\V\|_F=\sqrt{\V\ast_{d+1}\V}$ and $E_1\in\R^m$ is the 
first vector from the canonical basis, analogously to the approach used when $\A$ is a 
square matrix and $\V$ is a vector; see \cite{BR,GIJB,HYM}. 

This approach to approximate $f(\A)\ast_{d+1}\V$ works well when $\V=\E$ and can be applied to 
determine accurate approximations of the multilayer total communicability \eqref{mltc} and
the multilayer total communicability of node $v_i^{\ell_1,\ldots , \ell_d}$ defined by \eqref{mltciell}. The
former is approximated by 
\[
\E\ast_{d+1} \VV_m\ast_1f(H_m)\ast_1E_1 \|\V\|_F 
\]
and the latter by 
\begin{equation}\label{mltciellred}
\E_{i,\ell_1,\ldots , \ell_d}\ast_{d+1} \VV_m\ast_1f(H_m)\ast_1E_1 \|\V\|_F.
\end{equation}
In particular, the evaluation of \eqref{mltciellred} does not require any arithmetic work
when the expression \eqref{fred} is available. This makes the evaluation of the multilayer 
total communicability of all nodes $v_i^{\ell_1,\ldots , \ell_d}$, $1\leq i\leq N$ and $1\leq\ell_s\leq K_s,~~,1\leq s\leq d$, 
inexpensive when the expression \eqref{fred} is known. We use this fact when determining 
nodes for which this measure is large in Section \ref{sec5}. 

However, Algorithm \ref{Algo1} often suffers from breakdown when seeking
to approximate an expression of the form $f(\A)\ast_{d+1}\V$ when the tensor $\V$ is sparse, i.e.,
when $\V$ has many vanishing entries. This is the case when seeking to approximate the 
subgraph centrality \eqref{sgc} by
\[
\E_{i,\ell_1,\ldots , \ell_d}\ast_{d+1} \VV_m\ast_1 f(H_m)\ast_1 E_1 \| \E_{i,\ell_1,\ldots , \ell_d}\|_F,
\]
or the communicability \eqref{com} between the node $v_i^{\ell_1,\ldots , \ell_d}$ and node $v_j^{k_1,\ldots , k_d}$ by 
\[
\E_{i,\ell_1,\ldots , \ell_d}\ast_{d+1} \VV_m\ast_1 f(H_m) \ast_1 E_1 \| \E_{j,k_1,\ldots , k_d}\|_F.
\]
When computing these approximations, the initial block tensor is $\V=\E_{i,\ell_1,\ldots , \ell_d}$, which is 
very sparse. Since the tensor $\A$ typically also is sparse, this often results in that
the scalar $h_{j+1,j}$ in line 10 of Algorithm \ref{Algo1} vanishes for some 
$1\leq j\leq m$. The computations with the algorithm then cannot be continued, and the
available expression  at breakdown,
\[
\VV_j\ast_1f(H_j)\ast_1E_1 \|\V\|_F,
\] 
might not furnish an approximation of desired accuracy. Moreover, even when all the 
multilayer subgraph centralities can be computed to determine the node with the largest 
subgraph centrality, this is quite expensive for large multilayer networks. We describe 
in the following subsection how these difficulties can be reduced by replacing the 
initial tensor $\V$ in Algorithm \ref{Algo1} by a block of tensors. 

\subsection{A block Arnoldi process based on the Einstein product}
We describe a block Arnoldi process that uses the Einstein product. It differs from 
Algorithm \ref{Algo1} in that the initial tensor $\V$ is extended to a block tensor. The
application of $m$ steps of the block Arnoldi process to $\A$ with initial tensor $\W$
determines, in the absence of breakdown, an orthonormal basis for the block tensor Krylov 
subspace
\begin{align}
\K_{m+1}^{\rm block}(\A,\W)&={\rm range}\{\W,\A\ast_{d+1} \W,\ldots ,\A^m\ast_{d+1} \W\} \label{BlockKrSub}\\ 
&=\left\{\sum_{i=0}^m\A^i\ast_{d+1}\W\ast_1\Omega_i,~ \Omega_i\in\R^{P\times P}\right\},\nonumber
\end{align}
where $\W\in\R^{N\times K_1\times\ldots\times K_d\times P}$. We will refer to the integer $P$ as the block size.   
This definition of the subspace is analogous to the definition of the
solution subspace for block matrix methods used in \cite{JST}.
The block Arnoldi process of this subsection has the advantage of typically requiring 
fewer accesses to the adjacency tensor than when applying Algorithm 
\ref{Algo1}. Moreover, choosing suitable auxiliary columns in the initial tensor $\W$, 
the occurrences of breakdowns can be reduced in comparison with Algorithm \ref{Algo1}. 
The block Arnoldi process is summarized by Algorithm \ref{Algo2}. 

\begin{algorithm}
\caption{Block tensor Arnoldi process}\label{Algo2}
\begin{algo}
\INPUT Adjacency tensor $A\in\R^{{N\times K_1\times\ldots\times K_d\times N\times K_1\times\ldots\times K_d }}$, initial tensor 
$\W\in\R^{N\times K_1\times\ldots\times K_d\times P}$, and number of steps $m$.
\OUTPUT Orthonormal basis $\WW_{m+1}=\{\W_1,\W_2,\ldots,\W_{m+1}\}$ for the block
Krylov subspace \eqref{BlockKrSub}, and nontrivial entries $H_{i,j}\in\R^{P\times P}$ of
the upper block Hessenberg matrix $\HH_{m+1,m}=[H_{i,j}]\in\R^{P(m+1)\times Pm}$.
\STATE Compute the QR factorization $\W=\QQ\ast_1 R$, where the tensor 
$\QQ\in\R^{N\times K_1\times\ldots\times K_d\times P}$ satisfies $\QQ^T\ast_1 \QQ = \I$ and the matrix
$R\in\R^{P\times P}$ is upper triangular. Set $\W_1=\QQ$ and $H_{1,0}=R$.
\FOR $j=1,\ldots ,m$ do
\STATE $\U = \A\ast_{d+1}\W_j$
\FOR $i=1,\ldots ,j$ do 
\STATE $H_{i,j} = \W_i^T\ast_{d+1}\U,$ 
\STATE $\U = \U-\W_i\ast_1H_{i,j}$
\ENDFOR
\STATE Compute the QR factorization $\U=\QQ\ast_1 R$, where $\QQ\in\R^{N\times K_1\times\ldots\times K_d\times P}$
satisfies $\QQ^T\ast_1 \QQ = \I$ and the matrix $R\in\R^{P\times P}$ is upper triangular.
Set $\W_{j+1}=\QQ$ and $H_{j+1,j}=R$.
\ENDFOR
\end{algo}
\end{algorithm} 

Algorithm \ref{Algo2} determines an orthonormal basis 
\[
\WW_{m+1}=\left[\W_1,\ldots,\W_{m+1}\right]\in\R^{N\times K_1\times\ldots\times K_d\times P(m+1)}
\]
for the block Krylov subspace \eqref{BlockKrSub} and the upper block Hessenberg matrix 
\[
\HH_{m+1,m}=\left[\begin{array}{cccccc} 
H_{1,1} & H_{1,2} & H_{1,3} & \cdots & H_{1,m-1} & H_{1,m} \\
H_{2,1} & H_{2,2} & H_{2,3} & \cdots & H_{2,m-1} & H_{2,m} \\
        & H_{3,2} & H_{3,3} & \cdots & H_{3,m-1} & H_{3,m} \\
	&         & \ddots  &        & \vdots    & \vdots \\
        &         &         &        & H_{m,m-1} & H_{m,m} \\
        &         &         &        &           & H_{m+1,m} 
\end{array}\right]\in\R^{P(m+1)\times Pm}.
\]
Its leading $Pm\times Pm$ submatrix is denoted by $\HH_m$.
We have the following result:

\begin{pro}
Suppose that $m$ steps of Algorithm \ref{Algo2} have been carried out. Then
\begin{equation}\label{Block1}
\A\ast_{d+1}\WW_m=\WW_{m+1}\ast_1\HH_{m+1,m}
\end{equation}
and
\begin{equation}\label{Block2}
\A\ast_{d+1}\WW_m=\WW_m\ast_1\HH_m+\W_{m+1}\ast_1H_{m+1,m}\ast_1\EE_m^T,
\end{equation}
where $\WW_m=[\W_1,\ldots,\W_m]\in\R^{N\times K_1\times\ldots\times K_d\times Pm}$ is made up of the first $m$ tensor columns of $\WW_{m+1}$ and $\EE_m^T=[0,0,\ldots , I_P]\in\R^{P\times Pm}$. Moreover, 
\begin{equation}\label{Block3}
\A^p\ast_{d+1}\W = \WW_m\ast_1\HH_m^p\ast_1E_1\chi_0,
\end{equation}
where $\chi_0$ is obtained from the QR factorization of $\W$, such that 
$\W=\WW_m\ast_1E_1\chi_0$ holds for all $p\leq0$. 
\end{pro}
\noindent
Proof: 
For $1\leq j\leq m$, one has
\begin{align*}
&[\WW_{m+1}\ast_1\HH_{m+1,m}]_{:,\ldots ,:,1+P(j-1):Pj}& \\
~~~&= \sum_{i=1}^{j+1}(\WW_{m+1})_{:,\ldots ,:,1+P(i-1):Pi} \ast_1(\HH_{m+1,m})_{1+P(i-1):Pi,1+P(j-1):Pj}\\
~~~&=[\A\ast_{d+1}\WW_m]_{:,\ldots ,:,1+P(j-1):Pj}.
\end{align*}
This leads to equation \eqref{Block1}; equation \eqref{Block2} can be shown similarly.
The last claim can be proved by induction. In fact, if $p=0$, then one has 
$\A^0\ast_{d+1}\W = \W = \WW_m\ast_1E_1\chi_0$ and, by assuming that \eqref{Block3} holds for
$p\geq0$, one obtains
$$
\A^{p+1}\ast_{d+1}\W = \A\ast_{d+1}\A^p\ast_{d+1}\W = \A\ast_{d+1}\WW_m\ast_1\HH_m^p\ast_1E_1\chi_0,
$$
so that, using equation \eqref{Block2}, we have
\begin{eqnarray*}
\A^{p+1}\ast_{d+1}\W &=&(\WW_m\ast_1\HH_m+\W_{m+1}\ast_1H_{m+1,m}\ast_1\EE_m^T)
\ast_1\HH_m^p\ast_1E_1\chi_0\\
&=&\WW_m\ast_1\HH_m\ast_1\HH_m^p\ast_1E_1\chi_0\\
&&+\W_{m+1}\ast_1H_{m+1,m} \ast_1\EE_m^T\ast_1\HH_m^p\ast_1E_1\chi_0,
\end{eqnarray*}
where the second term vanishes due to the fact that $\HH_m$ is a block Hessenberg matrix. 
This concludes the proof. $~~~\Box$
 
Thanks to equation \eqref{Block3}, $f(\A)\ast_{d+1}\W$ can be approximated by 
\begin{equation}\label{apprx}
\WW_m\ast_{d+1}f(\HH_m)\ast_1E_1\chi_0,
\end{equation}
where $\chi_0$ is such that $\W=\WW_m\ast_1E_1\chi_0$. The advantage of this approach is 
that we can compute the multilayer subgraph centrality and the resolvent-based subgraph 
centrality of $P$ nodes at once, and we also can determine approximations of the 
multilayer communicabilities of these $P$ nodes essentially for free. This is because we
are approximating the quantity $f(\A)\ast_{d+1}\W$ by \eqref{apprx}, and then can evaluate 
the approximation $\W^T\ast_{d+1}\WW_m\ast_{d+1}f(\HH_m)\ast_1 E_1\chi_0$ of 
$\W^T\ast_{d+1}f(\A)\ast_{d+1}\W$ inexpensively. Moreover, we can circumvent the numerical 
stability issue related to the sparsity of the adjacency tensor and the initial block by 
adding a dense tensor in the initial block tensor. This technique has been discussed for 
matrix functions case in \cite{FMRR}. Notice that if we include $\E$, i.e., the tensor of 
all ones, in the initial block tensor, then Algorithm \ref{Algo2} will produce the same 
quantities of interest as one would compute with Algorithm \ref{Algo1} when applied with 
the initial tensor $\V=\E$, as well as the multilayer total network communicability. This 
is because once we approximate 
\[
f(\A)\ast_{d+1}\W:=[f(\A)\ast_{d+1}\W_1,\ldots,
f(\A)\ast_{d+1}\W_P,f(\A)\ast_{d+1}\W_{P+1}],
\]
where $\W_{P+1}=\E$, we only need to compute the following Einstein product 
$M:=\W^T\ast_{d+1}f(\A)\ast_{d+1}\W\in\R^{P+1\times P+1}$. Then for $1\leq i,j\leq P$ and 
$i\neq j$ the quantities $M_{i,j}$ are the communicabilities between different nodes, for
$1\leq i\leq P$ the quantities $M_{i,i}$ are the multilayer subgraph centralities, for 
$1\leq i\leq P$ the quantities $M_{i,P+1}$ are the same obtained by Algorithm \ref{Algo1},
and $M_{P+1,P+1}$ is the multilayer total network communicability of the whole network.

\section{Computed examples}\label{sec5}
This section presents some examples to illustrate the performance of the methods discussed
above. The computations were carried out using MATLAB R2015b. We use the Matlab library, tensor toolbox \cite{KB1}, to perform operations on tensors. For the examples in Sections
\ref{ex2} and \ref{ex3}, we choose the minimum number of steps, $m$, with the Krylov 
subspace method needed to obtain the same ranking of the first $10$ nodes as the ranking 
obtained when evaluating the exact tensor function. We will see that the number of steps
required is quite small. Due to the size of the networks in the examples in Sections 
\ref{ex4} and \ref{ex5}, it is expensive to evaluate the exact tensor function. In these
examples, we therefore increase the number of steps, $m$, of the Krylov subspace until 
the ranking does not change, and consider the ranking so obtained the exact one.
Knowledge of the nodes for which we are interest in computing multilayer subgraph 
centrality is needed for Algorithm \ref{Algo2}. We choose these nodes as the top central 
ones obtained by Algorithm \ref{Algo1}.

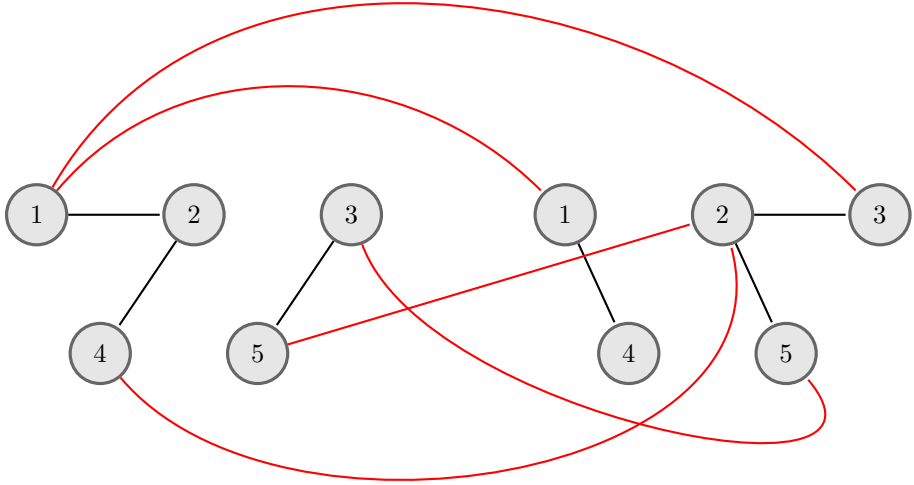
\begin{figure}[H]
\begin{center}
\begin{tikzpicture}[
 roundnode/.style={circle, draw=black!60, fill=black!10, very thick, minimum size=8mm},
 ->,>=stealth',shorten >=1pt,node distance=2.5cm,auto,thick]
\node[roundnode]      (node1)                                                        {$1$};
\node[roundnode]      (node2)       [right=1.25cm of node1]           {$2$};
\node[roundnode]      (node3)       [right=1.25cm of node2]           {$3$};
\node[roundnode]      (node4)       [below right=1.25cm and 0.25cm of node1]         {$4$};
\node[roundnode]      (node5)       [right=1.25cm of node4]            {$5$};

\node[roundnode]      (node6)        [right=2cm of node3]                                                  {$1$};
\node[roundnode]      (node7)       [right=1.25cm of node6]           {$2$};
\node[roundnode]      (node8)       [right=1.25cm of node7]          {$3$};
\node[roundnode]      (node9)       [below right=1.25cm and 0.25cm of node6]         {$4$};
\node[roundnode]      (node10)       [right=1.25cm of node9]            {$5$};

\draw[-] (node1) -- (node2);
\draw[-] (node2) -- (node4);
\draw[-] (node3) -- (node5);

\draw[-] (node6) -- (node9);
\draw[-] (node7) -- (node8);
\draw[-] (node7) -- (node10);

\rd{
\draw[-] (node1) to [out=50] (node6);
}
\rd{
\draw[-] (node1) to [out=60] (node8);
}
\rd{
\draw[-] (node3) to [out=-70,in=-50] (node10);
}
\rd{
\draw[-] (node4) to [out=-50,in=-75] (node7);
}
\rd{
\draw[-] (node5) -- (node7);
}
\end{tikzpicture}
\end{center}
\caption{Example 1: Layers are presented from left to right in the order $L=1$ and $L=2$. 
The edges connecting nodes from same layer are marked in black. The edges 
connecting nodes from different layers are marked in red.}\label{fig1}
\end{figure}

\subsection{Example 1: A small synthetic multilayer network}\label{ex1}
Consider a small synthetic unweighted and undirected multilayer network with $d=1$, 
$K_1=2$, and $N=5$. It consists of $10$ edges as shown by the graph 
of Figure \ref{fig1}. Let $\A$ be the adjacency tensor. We compute the different 
multilayer centrality measures discussed using the Einstein product. The results are 
summarized in Table \ref{tab:table1} for $\beta=1$ and $\alpha=0.5/\lambda_{max}$. We can see 
that node $v_2^2$ has the largest centrality measures and, therefore, is 
the most important node. This node is followed by node $v_1^1$. Node $v_4^2$ is the
least important node. All centrality measures perform as can be expected. 
 
\begin{table}
\caption{Multilayer total communicability of the nodes $v_i^\ell$ (MTC), multilayer Katz
centrality (MKC), multilayer subgraph centrality with the modified tensor exponential 
(MSC$_{\exp_{0}}$),  multilayer subgraph centrality with the modified tensor resolvent 
(MSC$_{\rm res_0}$) for the nodes of Example 1.}
\label{tab:table1}
\centering
\begin{tabular}{c|cc|cc}
 & $MTC(i,\ell)$ & & $MKC(i,\ell)$ &\\
   \hline
 $v_i^{\ell}$ & $\ell=1$ & $\ell=2$  & $\ell=1$ & $\ell=2$    \\
    \hline
 $i=1$ & $12.2520$ & $7.7379$  & $2.1131$ & $1.7044$  \\    

 $i=2$ & $9.6537$ & $17.2450$  & $1.8145$ & $2.5454$  \\  

 $i=3$ & $8.9474$ & $11.7351$  & $1.7645$ & $1.9484$   \\  

 $i=4$ & $10.8250$ & $4.2550$  & $1.8876$ & $1.3470$   \\  

 $i=5$ & $10.6175$ & $10.6175$  & $1.8774$ & $1.8774$   \\  
    \hline
 &  $MSC_{\exp_{0}}(i,\ell)$ & & $MSC_{\rm res_0}(i,\ell)$ \\
   \hline
 $v_i^{\ell}$  & $\ell=1$ & $\ell=2$  & $\ell=1$ & $\ell=2$  \\
    \hline
 $i=1$ & $3.1001$ & $2.2834$ & $1.1507$ & $1.0952$  \\    

 $i=2$ & $2.3582$ & $4.1313$ & $1.0987$ & $1.2165$  \\  

 $i=3$ & $2.3946$ & $2.4698$ & $1.1003$ & $1.1039$  \\  

 $i=4$ & $2.4174$ & $1.5922$ & $1.1016$ & $1.0454$  \\  

 $i=5$ & $2.5001$ & $2.5001$ & $1.1051$ & $1.1051$  \\  
    \hline
\end{tabular}
\end{table}
 
\subsection{Example 2: A synthetic multilayer network}\label{ex2}
\begin{figure}
\centering
\begin{minipage}{.5\textwidth}
  \centering
  \includegraphics[width=\textwidth]{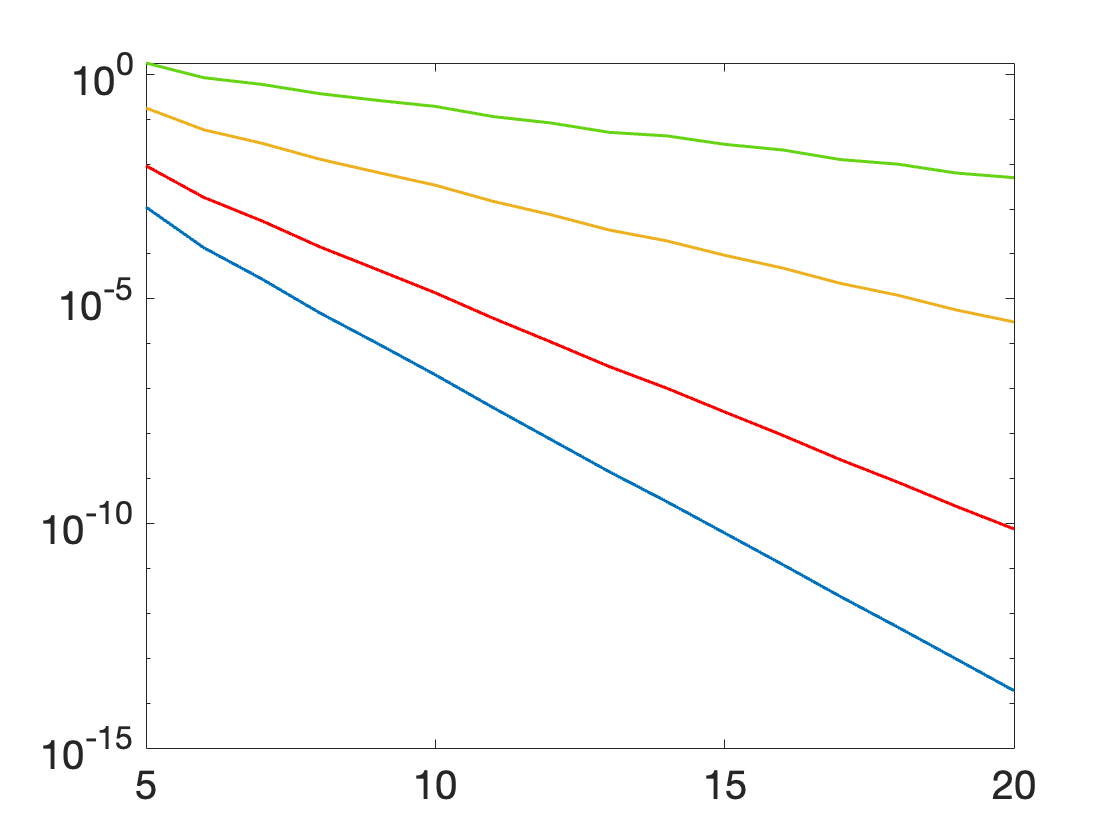}\\(a)
\end{minipage}%
\begin{minipage}{.5\textwidth}
  \centering
  \includegraphics[width=\linewidth]{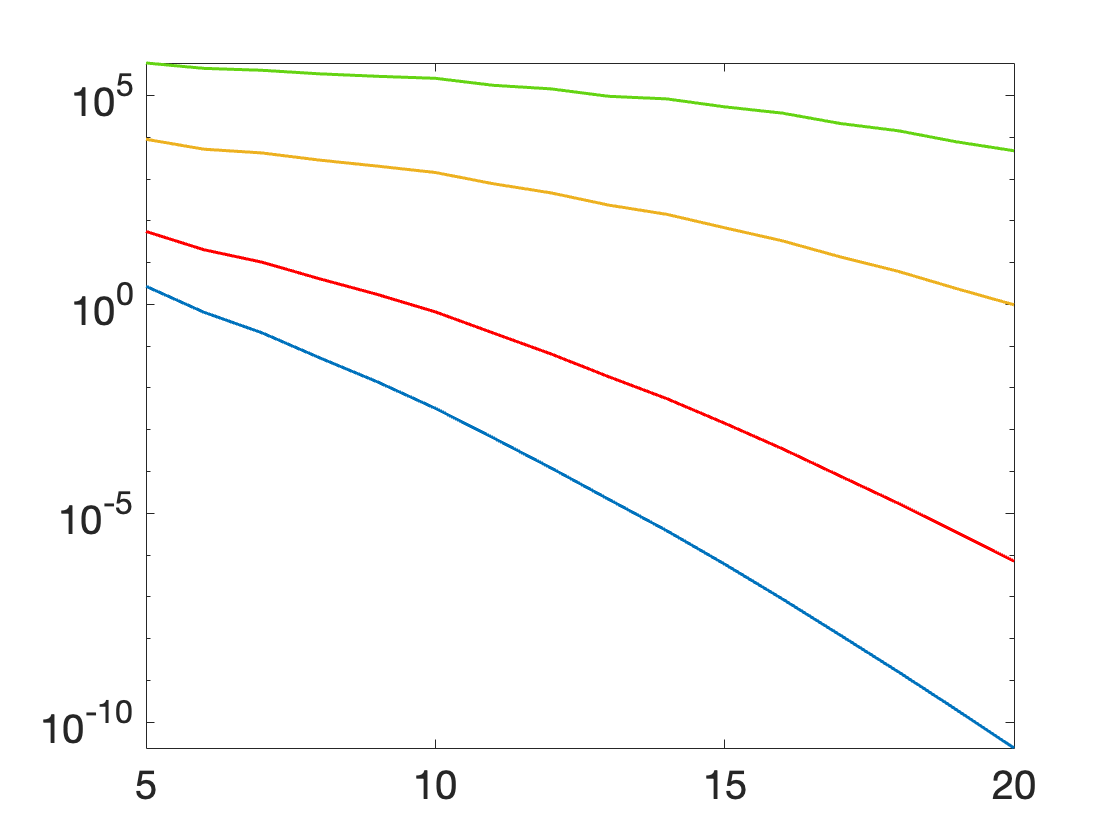}\\(b)
\end{minipage}
\newline
\centering
\begin{minipage}{.5\textwidth}
  \centering
  \includegraphics[width=\linewidth]{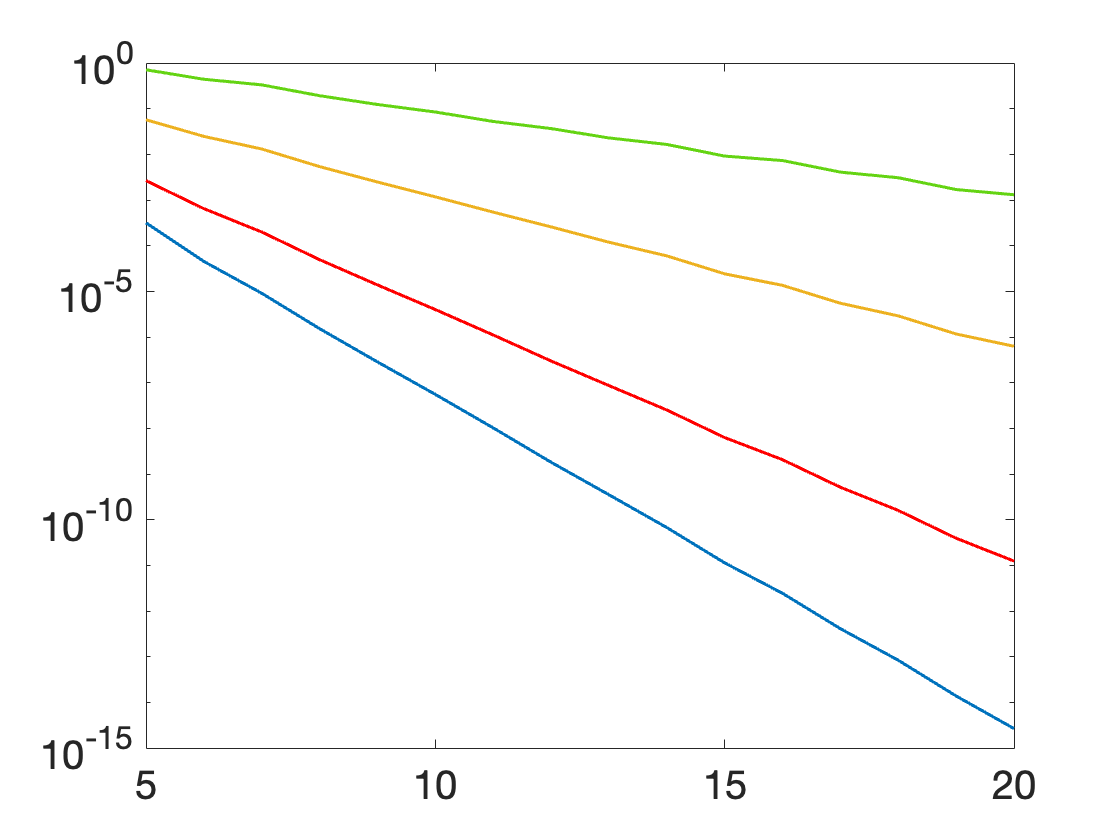}\\(c)
\end{minipage}%
\begin{minipage}{.5\textwidth}
  \centering
  \includegraphics[width=\linewidth]{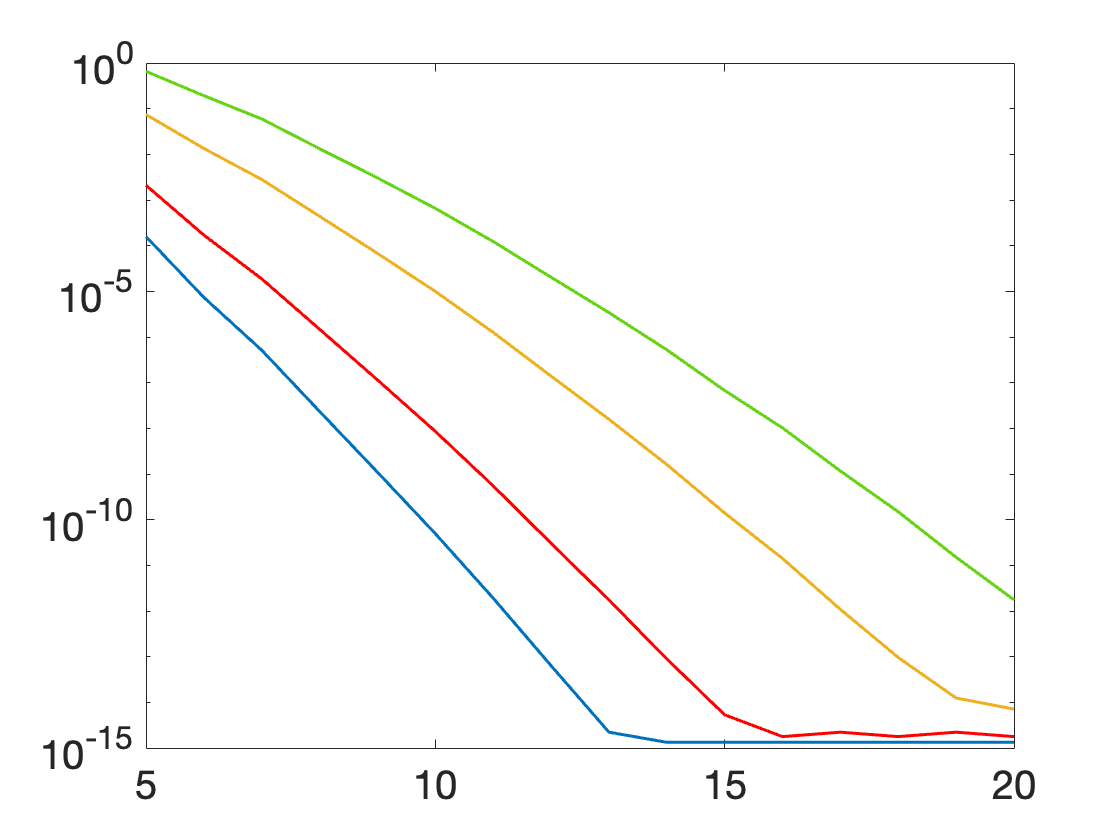}\\(d)
\end{minipage}
\caption{Infinity norm error for multilayer total communicability of the nodes and 
multilayer Katz centrality as functions of the Krylov subspace dimension $m$ for Example 
$2$ with $\alpha\in\{\frac{0.2}{\lambda_{\max}},
\frac{0.3}{\lambda_{\max}},\frac{0.5}{\lambda_{\max}}, \frac{0.7}{\lambda_{\max}}\}$, and
$\beta\in\{0.3,0.5,1,1.5\}$ ($\{$blue, red, yellow, green$\}$). 
$(a)$ multilayer Katz centrality for weighted edges, $(b)$ multilayer total 
communicability of the nodes for weighted edges, $(c)$ multilayer Katz centrality for 
unweighted edges, and $(d)$ multilayer total communicability of the nodes for unweighted 
edges.}
\label{fig2}
\end{figure}

We consider an example of a weighted multilayer network with $d=1$ aspect, $K_1=32$ 
layers, $N=20$ nodes, and $674$ edges. The network can be downloaded from  
\texttt{https://github.com/wjj0301/Multiplex-Networks}.\ We compute approximations of the
multilayer total communicability of the nodes $v_i^\ell$ and multilayer Katz centrality 
using Algorithm \ref{Algo1}, as well as the exact multilayer total communicability of the 
nodes $v_i^\ell$ and multilayer Katz centrality, i.e., the subgraph centrality determined
with the function \eqref{modres}, by flattening the adjacency tensor and 
using the MATLAB function {\sf expm} and the MATLAB backslash operator when using the
function \eqref{modres}. The errors in the 
computed approximations are determined by the vector infinity norm, that is, the exact 
multilayer total communicability of the nodes, multilayer Katz centrality, and their 
approximations are stored in vectors and the infinity norm is applied to measure the 
distance between the vectors with the exact entries and the vectors with the corresponding
entries determined by Krylov subspace methods. Figure \ref{fig2} displays the errors as a 
function of the Krylov subspace dimension $m$. Table \ref{tab:table2} lists the top $10$ 
central nodes according to the multilayer total communicability with respect to a node and
multilayer Katz centrality. We notice that the multilayer Katz centrality is approximated 
accurately when the same number of steps, $m$, of Algorithm \ref{Algo1} are carried out. 
The multilayer total communicability of the nodes is approximated less accurately with the
same number of steps, but still gives the same ranking of nodes as when the exact 
communicability is used.

 \begin{table}
\caption{Top $10$ central nodes determined by Algorithm \ref{Algo1} with $m=6$ using the 
(approximate) multilayer total communicability of the nodes (MTC) and (approximate) 
multilayer Katz centrality (MKC) so obtained for the weighted multilayer network in 
Example 2, with $\alpha=0.4/\lambda_{\max}$ and  $\beta=0.4$. These computed values are 
compared to the exact ones obtained by flattening and using matrix functions.}
\label{tab:table2}
\centering
{\small
\begin{tabular}{ccc|ccc}
$\{i,\ell\}$ &  $MTC\{i,\ell\}$ & $MTC_{\rm exact}\{i,\ell\}$ & $\{i,\ell\}$ & $MKC\{i,\ell\}$ 
& $MKC_{\rm exact}\{i,\ell\}$ \\
   \hline
$\{18,24\}$ & $99.4969$   & $97.1435$ & $\{18,24\}$ & $3.6474$   & $3.6407$   \\

$\{17,26\}$ & $70.0610$   & $72.6593$ & $\{19,19\}$ & $3.0148$  & $3.0191$   \\

$\{13,19\}$ & $64.5009$   & $67.7357$ & $\{14,23\}$ & $2.8966$   & $2.8956$   \\
    
$\{8,26\}$ & $62.8084$   & $67.1222$ & $\{13,19\}$ & $2.6318$   & $2.6407$   \\

$\{19,19\}$ & $62.1558$  & $63.7426$ & $\{17,26\}$ & $2.6242$   & $2.6313$   \\
  
$\{6,24\}$ & $59.3319$   & $60.2469$ & $\{6,29\}$ & $2.5999$   & $2.5990$   \\

$\{19,4\}$ & $58.9334$   & $59.4300$ & $\{2,29\}$ & $2.5497$   & $2.5477$   \\

$\{14,24\}$ & $56.5841$   & $56.4157$ & $\{1,32\}$ & $2.5371$   & $2.5425$   \\

$\{1,32\}$ & $52.8265$   & $54.7537$ & $\{1,3\}$ & $2.5251$   & $2.5240$   \\

$\{2,24\}$ & $51.6585$  & $54.3109$ & $\{14,24\}$ & $2.5150$   & $2.5146$   \\
     \hline
\end{tabular}
}
\end{table}

We also set the weight of all edges to one to obtain an unweighted multilayer network 
and compute the multilayer total communicability of the nodes and multilayer Katz 
centrality. We notice that the multilayer total communicability of the nodes is approximated 
accurately when $m=6$ steps of Algorithm \ref{Algo1} are carried out also for larger 
values of $\beta$. However, the multilayer Katz centrality is not approximated accurately 
for $m=6$ steps when $\alpha>0.66/\lambda_{\max}$. Figure \ref{fig2} depicts the influence
of the values of $\alpha$ and $\beta$ on the dimension of Krylov subspace needed to obtain 
accurate approximations for both weighted and unweighted graphs.

Finally, we turn to the computation of multilayer subgraph centralities of some nodes. We 
apply Algorithm \ref{Algo2} to our adjacency tensor 
with the first block $\V\in\R^{N\times K_1\times P}$ with $P$ frontal slices of the form 
$E_{i,\ell}$. We let $P=10$ and choose the indices $\{i,\ell\}$ randomly. The results 
obtained are summarized in Table \ref{tab:table3}, which also shows communicabilities 
that we obtain for free since we use a block method.

\begin{table}
\caption{Multilayer subgraph centrality obtained with the modified tensor exponential 
(MSC$_{\exp_{0}}$) and multilayer subgraph centrality with the modified tensor resolvent 
(MSC$_{\rm res_0}$) for some nodes for the unweighted multilayer network in Example 2, as
well as free multilayer communicabilities of the nodes (MC) with 
$\alpha=0.5/\lambda_{\max}$, $\beta=1$, and $m=5$.}\label{tab:table3}
\centering
\begin{tabular}{ccc|cc}
   \hline
$\{i,\ell\}$ & $MSC_{\exp_{0}}\{i,\ell\}$ & $MSC_{\rm res_0}\{i,\ell\}$ & $\{i,\ell,j,k\} $& 
$MC\{i,\ell,j,k\} $ \\
     \hline
  $\{4,1\}$ & $1.1683$ & $1.0358$ & & \\
  $\{5,1\}$ & $1.2116$ & $1.0488$ & & \\
  $\{12,1\}$ & $1.5878$ & $1.1321$ & &\\
  $\{1,3\}$ & $1.1766$ & $1.0395$ & $\{5,1,4,1\}$ & $1.0505$ \\
  $\{19,4\}$ & $1.0084$ & $1.0038$ &  $\{12,1,4,1\}$ & $0.1768$  \\
  $\{2,11\}$ & $1.0592$ & $1.0193$ &  $\{12,1,5,1\}$ & $0.5531$  \\
  $\{13,18\}$ & $1.1681$ & $1.0354$ & $\{13,18,19,18\}$ & $0.0502$  \\
  $\{19,18\}$ & $1.0097$ & $1.0048$ &  & \\
  $\{14,22\}$ & $1.2537$ & $1.0611$ & & \\
  $\{8,25\}$ & $1.5431$ & $1.1176$ &  & \\
     \hline      
     
\end{tabular}
\end{table}

\subsection{Example 3: Multiplex network (Scotland Yard transportation data)}\label{ex3}
\begin{figure}
\centering
\begin{minipage}{.5\textwidth}
  \centering
\includegraphics[width=\linewidth]{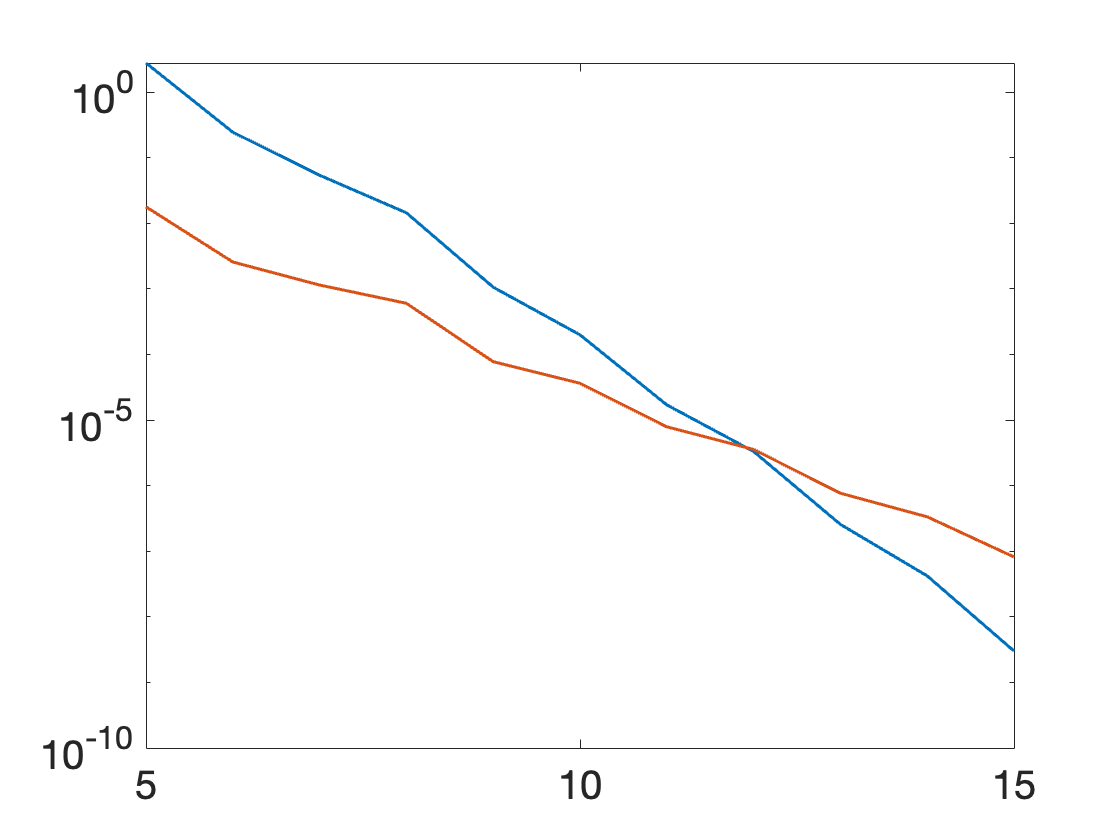}\\(a)
\end{minipage}%
\begin{minipage}{.5\textwidth}
  \centering
\includegraphics[width=\linewidth]{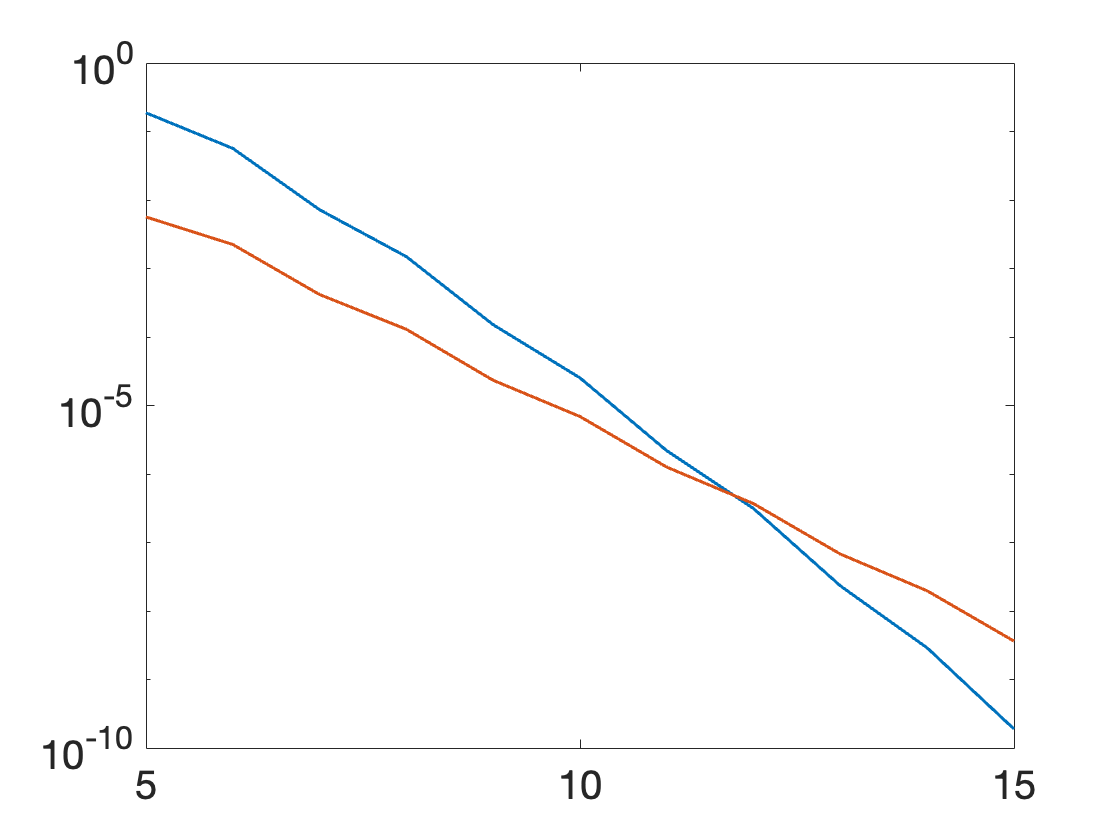}\\(b)
\end{minipage}
\caption{Infinity norm error for multilayer total communicability of the nodes (blue) and 
multilayer Katz centrality (orange) as functions of the Krylov subspace dimension $m$ for 
Example $3$. $(a)$ when edges are unweighted and $\alpha=0.4/\lambda_{\max}$ and 
$\beta=0.4$, $(b)$ when edges are weighted and $\alpha=0.3/\lambda_{\max}$ 
and $\beta=0.1$.}
\label{fig3}
\end{figure}

This example considers the Scotland Yard transportation network created by the authors 
of \cite{BS}, which is a multiplex network. A multiplex network is a special case of a 
multilayer network. The network can be downloaded from \cite{B_repository} as a weighted 
or unweighted multiplex network. It consists of $3324$ edges and $N=199$ nodes that
represent public transport stops in the city of London. The network has  $K_1=4$ layers that 
represent different modes of transportation: Boat, underground, bus, and taxi. The weights 
are determined so that the edges in the layer that represent travel by taxi all have 
weight one. A taxi ride is defined as a trip by taxi between two adjacent nodes in the 
taxi layer; a taxi ride along $k$ edges is considered $k$ taxi rides. The weights of edges
in the boat, underground, and bus layers are chosen to be equal to the minimal number of 
taxi rides required to travel between the same nodes. We compute the multilayer Katz 
centrality and multilayer total communicability of the nodes by applying Algorithm 
\ref{Algo1} to the adjacency tensor of the given multiplex network, as well as to the 
adjacency tensor of the associated unweighted network; in the latter adjacency tensor all 
edges have weight one. We compute the exact  multilayer Katz centrality and multilayer 
total communicability of the nodes in the same way as in the previous example and evaluate
the infinity norm error as well. Figure \ref{fig3} displays the errors as a function of 
the Krylov subspace dimension. This illustrates the accuracy of Algorithm \ref{Algo1} when 
applied to multiplex networks.

We also compute the multilayer subgraph centrality (MSC) measures for the top $10$ nodes using Algorithm 
\ref{Algo2}; see Table \ref{tab:table4}. The total network communicability obtained is 
$3.904\cdot 10^3$. Multilayer communicabilities between some of the nodes considered are listed in Table \ref{tab:table5}.

\begin{table}
\caption{Multilayer subgraph centrality obtained with the modified tensor exponential 
(MSC$_{\exp_{0}}$) and multilayer subgraph centrality with the modified tensor resolvent 
(MSC$_{\rm res_0}$) for some nodes of the network in Example 3, with 
$\alpha=0.3/\lambda_{\max}$, $\beta=0.3$, $P=11$ and $m=10$.}

\label{tab:table4}
\centering
{\small
\begin{tabular}{ccc}
$\{i,\ell\}$ &  $MSC_{\exp_{0}}\{i,\ell\}$ & $MSC_{\rm res_0}\{i,\ell\}$  \\
   \hline
$\{\text{142,4}\}$ & $1.6080$  &  $1.0345$     \\

$\{\text{140,4}\}$ & $1.5004$   &   $1.0295$     \\

$\{\text{58,4}\}$ & $1.5182$   &   $1.0301$     \\
    
$\{\text{128,4}\}$ & $1.4583$   &   $1.0267$     \\

$\{\text{67,3}\}$ & $1.4442$  &   $1.0262$     \\

$\{\text{153,4}\}$ & $1.4457$   &   $1.0263$     \\

$\{\text{143,4}\}$ & $1.4754$   &   $1.0273$      \\

$\{\text{114,4}\}$ & $1.4933$   &   $1.0293$     \\

$\{\text{129,4}\}$ & $1.4647$   &   $1.0269$     \\

$\{\text{128,3}\}$ & $1.4519$  &   $1.0265$     \\
     \hline
\end{tabular}
}
\end{table}

\begin{table}
\caption{Multilayer communicabilities obtained with the modified tensor exponential 
(MC$_{\exp_{0}}$) and with the modified tensor resolvent 
(MC$_{\rm res_0}$) for some nodes of the network in Example 3, with 
$\alpha=0.3/\lambda_{\max}$, $\beta=0.3$, $P=11$ and $m=10$.}

\label{tab:table5}
\centering
{\small
\begin{tabular}{ccc}
$\{i,\ell,j,k\}$ &  $MC_{\exp_{0}}$ & $MC_{\rm res_0}$  \\
   \hline
$\{\text{128,4,142,4}\}$ & $0.4666$  &  $0.0614$     \\

$\{\text{153,4,140,4}\}$ & $0.1267$  &  $0.0070$     \\

$\{\text{142,4,143,4}\}$ & $0.5192$  &  $0.0646$     \\

$\{\text{128,3,128,4}\}$ & $0.5010$   &   $0.0640$      \\

$\{\text{67,3,153,4}\}$ & $7.0328\cdot10^{-4}$   &   $1.8124\cdot 10^{-4}$     \\

   \hline
\end{tabular}
}
\end{table}

\subsection{Example 4: Multiplex network (European airlines data set)}\label{ex4}
The European airlines data set consists of $N=450$ nodes that represent European airports 
and has $K_1=37$ layers that represent different airlines operating in Europe. There are 
$3588$ edges, which represent available routes. This network can be represented by a 
fourth-order adjacency tensor $\A\in\R^{N\times K_1\times N\times K_1}$ such that 
$\A(i,\ell ,j,\ell)=1$ if there is a flight connecting airports $i$ and $j$ with airline 
$\ell$. Moreover, $\A(i,\ell,i,k)=1$ for every $1\leq \ell,k\leq K_1$ to reflect the effort 
required to change airlines for connecting flights. The network can be downloaded 
from \cite{B_repository}. Similarly as Taylor et al. \cite{TPM}, we only include $N=417$ 
nodes from the largest connected component of the network. We compute the multilayer total
communicability of the nodes and the multilayer Katz centrality using Algorithm 
\ref{Algo1} to approximate the tensor exponential and the tensor resolvent functions. 
Table \ref{tab:table6} lists the top $10$ central nodes. We obtained similar ranking as 
reported in \cite{BS}, where the authors applied Krylov subspace methods to the 
supra-adjacency matrix of the network in order to compute matrix function-based centrality
measures such as the Katz centrality.

We apply Algorithm \ref{Algo2} to our adjacency tensor in order to compute the 
multilayer subgraph centrality (MSC) for the top $10$ nodes determined earlier. This 
algorithm determines the multilayer subgraph centrality for the different nodes at once, the results are reported in Table \ref{tab:table7}. Table \ref{tab:table8} shows the multilayer communicabilities between the considered nodes
determined by the algorithm. The multilayer total network communicability is $2.4163.10^{7}$.

\begin{table}
\caption{Top $10$ central nodes according to multilayer total communicability of the nodes
(MTC) and multilayer Katz centrality (MKC) for the European airlines network in Example 4, 
with $\alpha=0.5/\lambda_{\max}$ and $\beta=0.2$. The MTC and MKC values are computed with 
Algorithm \ref{Algo1} with $m=20$. }
\label{tab:table6}
\centering
{\small
\begin{tabular}{cc|cc}
$\{i,\ell\}$ &  $MTC\{i,\ell\}$  & $\{i,\ell\}$ & $MKC\{i,\ell\}$  \\
   \hline
$\{\text{Stansted, Ryanair}\}$ & $8.2164\cdot 10^3$   &  $\{\text{Stansted, Ryanair}\}$ 
& $4.4228$     \\

    $\{\text{Munich, Lufthansa}\}$ & $7.5209\cdot 10^3$   &  $\{\text{Munich, Lufthansa}\}$ 
    & $4.0937$     \\

$\{\text{Frankfurt, Lufthansa}\}$ & $7.4541\cdot 10^3$   &  $\{\text{Frankfurt, Lufthansa}\}$ 
& $4.0650$     \\
    
$\{\text{Dublin, Ryanair}\}$ & $6.8942\cdot 10^3$   &  $\{\text{Ataturk, Turkish}\}$ 
& $4.0486$     \\

$\{\text{Gatwick, EasyJet}\}$ & $6.5714\cdot 10^3$  &  $\{\text{Gatwick, EasyJet}\}$ 
& $3.7925$     \\

$\{\text{Ataturk, Turkish}\}$ & $6.4399\cdot 10^3$   &  $\{\text{Dublin, Ryanair}\}$ 
& $3.6479$     \\

$\{\text{Amsterdam, KLM}\}$ & $6.0994\cdot 10^3$   &  $\{\text{Vienna, Austrian}\}$ 
& $3.5939$      \\

$\{\text{Vienna, Austrian}\}$ & $5.8057\cdot 10^3$   &  $\{\text{Amsterdam, KLM}\}$ 
& $3.5661$     \\

$\{\text{Caravaggio, Ryanair}\}$ & $5.7104\cdot 10^3$   &  $\{\text{Caravaggio, Ryanair}\}$ 
& $3.3244$     \\

$\{\text{Adolfo, Ryanair}\}$ & $5.5765\cdot 10^3$  &  $\{\text{Charles de Gaulle, Air France}\}$ 
& $3.2444$     \\
     \hline
\end{tabular}
}
\end{table}

\begin{table}
\caption{Multilayer subgraph centrality obtained with the modified tensor exponential 
(MSC$_{\exp_{0}}$) and multilayer subgraph centrality with the modified tensor resolvent 
(MSC$_{\rm res_0}$) for some nodes of the European airlines network in Example 4, with 
$\alpha=0.5/\lambda_{\max}$, $\beta=0.2$, $P=11$, and $m=10$.}

\label{tab:table7}
\centering
{\small
\begin{tabular}{ccc}
$\{i,\ell\}$ &  $MSC_{\exp_{0}}\{i,\ell\}$ & $MSC_{\rm res_0}\{i,\ell\}$  \\
   \hline
$\{\text{Stansted, Ryanair}\}$ & $53.6833$  &  $1.0282$     \\

$\{\text{Munich, Lufthansa}\}$ & $49.9200$   &   $1.0261$     \\

$\{\text{Frankfurt, Lufthansa}\}$ & $49.7844$   &   $1.0259$     \\
    
$\{\text{Dublin, Ryanair}\}$ & $49.0779$   &   $1.0225$     \\

$\{\text{Gatwick, EasyJet}\}$ &47.2554  &   $1.0239$     \\

$\{\text{Ataturk, Turkish}\}$ & $47.4633$   &   $1.0261$     \\

$\{\text{Amsterdam, KLM}\}$ & $45.1513$   &   $1.0224$      \\

$\{\text{Vienna, Austrian}\}$ & $45.2302$   &   $1.0228$     \\

$\{\text{Caravaggio, Ryanair}\}$ & $46.4080$   &   $1.0203$     \\

$\{\text{Adolfo, Ryanair}\}$ & $43.4764$  &   $1.0171$     \\
     \hline
\end{tabular}
}
\end{table}

\begin{table}
\caption{Multilayer communicabilities obtained with the modified tensor exponential 
(MC$_{\exp_{0}}$) and with the modified tensor resolvent 
(MC$_{\rm res_0}$) for some nodes of the European airlines network in Example 4, with 
$\alpha=0.5/\lambda_{\max}$, $\beta=0.2$, $P=11$, and $m=10$.}

\label{tab:table8}
\centering
{\small
\begin{tabular}{ccc}
$\{i,\ell,j,k\}$ &  $MC_{\exp_{0}}$ & $MC_{\rm res_0}$  \\
   \hline
$\{\text{Dublin, Ryanair, Stansted, Ryanair}\}$ & $1.33\cdot 10^1\phantom{-}$  & 
$2.17\cdot 10^{-2}$     \\

$\{\text{Vienna, Austrian, Stansted, Ryanair,}\}$ & $6.56\cdot 10^{-1}$  &  $4.37\cdot 10^{-5}$     \\

$\{\text{Frankfurt, Lufthansa, Munich, Lufthansa}\}$ & $1.31\cdot 10^{1\phantom{-}}$  &  
$2.46\cdot 10^{-1}$     \\

$\{\text{Amsterdam, KLM, Frankfurt, Lufthansa}\}$ & $5.58\cdot 10^{0\phantom{-}}$   &  
$8.96\cdot 10^{-4}$      \\

$\{\text{Caravaggio, Ryanair, Dublin, Ryanair}\}$ & $1.01\cdot 10^{1\phantom{-}}$   & 
$1.87\cdot 10^{-2}$     \\

$\{\text{Stansted, Ryanair, Ataturk, Turkish}\}$ & $7.42\cdot 10^{-1}$  &  $4.88\cdot 10^{-5}$     \\

$\{\text{Munich, Lufthansa, Gatwick, EasyJet}\}$ & $3.53\cdot 10^{0\phantom{-}}$   &  
$2.24\cdot 10^{-4}$     \\

$\{\text{Caravaggio, Ryanair, Adolfo, Ryanair}\}$ & $7.94\cdot 10^{0\phantom{-}}$   &  
$1.72\cdot 10^{-2}$     \\

$\{\text{Ataturk, Turkish, Caravaggio, Ryanair}\}$ & $5.06\cdot 10^{-1}$   &  
$3.20\cdot 10^{-5}$   \\  

$\{\text{Stansted, Ryanair, Adolfo, Ryanair}\}$ & $9.41\cdot 10^{0\phantom{-}}$   & 
$1.82\cdot 10^{-2}$   \\  
   \hline
\end{tabular}
}
\end{table}

\subsection{Example $5$: Wikispeedia network}\label{ex5}
This data set contains human navigation paths in Wikipedia, collected through the 
human-computation game Wikispeedia. Wikispeedia users are asked to navigate from a given 
source to a given target article, being allowed to click on links only. Nodes are articles of 
the English Wikipedia and edges represent clicks. The data is provided by the Stanford 
Network Analysis Project \texttt{http://snap.stanford.edu/index.html.}\ It contains $4604$ 
articles and $119882$ links. We classified the articles into $16$ different subjects 
(Countries, Science, Geography,~\ldots) and then built a multilayer network with $N=4604$
nodes and $K_1=16$ layers. Each layer contains the edges connecting the nodes that are considered 
to be classified in this layer. We have used the files that contain node identifiers and all 
edges from  \texttt{https://github.com/franloza/Wikispeedia-Network.}\ We compute the multilayer 
Katz centrality and multilayer total communicability of the nodes by applying Algorithm 
\ref{Algo1} to the adjacency tensor of the given multilayer network. Table 
\ref{tab:table9} lists the top $10$ central nodes for this network.

 \begin{table}
\caption{Top $10$ central nodes according to the multilayer total communicability (MTC) and 
the multilayer Katz centrality (MKC) for the Wikispeedia network in Example $5$, with $\beta=0.2$,
$\alpha=0.5/\lambda_{\max}$. The MTC and MKC values are approximations determined by 
Algorithm \ref{Algo1} with $m=20$.}
\label{tab:table9}
\centering
\begin{tabular}{cc}
$\{i,\ell\}$ &  $MTC\{i,\ell\}$  \\
   \hline
$\{\text{Fauna of Australia, Science}\}$ & $2.3711\cdot 10^{27}$  \\

    $\{\text{Africa, Countries}\}$ & $2.1895\cdot 10^{27}$  \\

$\{\text{Periodic table, Science}\}$ & $1.7854\cdot 10^{27}$   \\
   
$\{\text{List of elements by name, Science}\}$ & $1.7241\cdot 10^{27}$  \\
  
$\{\text{Periodic table (large version), Science}\}$ & $1.7162\cdot 10^{27}$ \\
 
$\{\text{Bird, Science}\}$ & $1.5887\cdot 10^{27}$  \\

$\{\text{Bird, Animals}\}$ & $1.5887\cdot 10^{27}$  \\

$\{\text{President of the United States, Politics}\}$ & $1.5476\cdot 10^{27}$  \\

$\{\text{Star, Science}\}$ & $1.5465\cdot 10^{27}$  \\

$\{\text{Star, Space}\}$ & $1.5465\cdot 10^{27}$   \\
     \hline
$\{i,\ell\}$ & $MKC\{i,\ell\}$  \\
   \hline
  $\{\text{Lebanon, Geography}\}$ & $4.7019$     \\

 $\{\text{Armenia, Countries}\}$ & $4.5967$     \\
 
 $\{\text{Armenia, Geography}\}$ & $4.5967$     \\
   
  $\{\text{Georgia, Countries}\}$ & $4.5712$     \\
  
 $\{\text{Georgia, Geography}\}$ & $4.5712$     \\

 $\{\text{Turkey, Countries}\}$ & $4.4296$     \\

  $\{\text{Turkey, Geography}\}$ & $4.4296$    \\

  $\{\text{Djibouti, Countries}\}$ & $4.2680$     \\

 $\{\text{Djibouti, Geography}\}$ & $4.2680$     \\

 $\{\text{Mozambique, Countries}\}$ & $4.1781$     \\
     \hline
\end{tabular}
\end{table}

\subsection{Example $6$: Synthetic multilayer network with $2$ aspects}\label{ex6}
In this last example, we consider a weighted undirected multilayer network with $d=2$ aspects, $K_1=3$, $K_2=2$ and $N=180$ nodes. The data can be download from \texttt{https://github.com/wjj0301/Multiplex-Networks},\  in the form of $6$ layers then, transformed to a multilayer network with $2$ aspects, it contains $148$ edges in total. We compute the multilayer Katz centrality and multilayer total communicability of the nodes by applying Algorithm \ref{Algo1} to the $6^{th}$ order adjacency tensor $\A\in\R^{180\times 3\times 2\times 180\times 3\times 2}$. Table \ref{tab:table10} displays the top $10$ central nodes. We apply Algorithm \ref{Algo2} to approximate the multilayer subgraph centrality measures for the top $10$ central nodes determined by Algorithm \ref{Algo1} and the total network communicability which is $1.6088\cdot 10^3$. Table \ref{tab:table11} displays multilayer subgraph centrality measures obtained by Algorithm \ref{Algo2}.

\begin{table}
\caption{Top $10$ central nodes according to multilayer total communicability of the nodes
(MTC) and multilayer Katz centrality (MKC) for the network in Example 6, 
with $\alpha=0.5/\lambda_{\max}$ and $\beta=0.2$. The MTC and MKC values are computed with 
Algorithm \ref{Algo1} with $m=10$.}
\label{tab:table10}
\centering
{\small
\begin{tabular}{cc|cc}
$\{i,\ell_1,\ell_2\}$ &  $MTC\{i,\ell_1,\ell_2\}$  & $\{i,\ell_1,\ell_2\}$ & $MKC\{i,\ell_1,\ell_2\}$  \\
   \hline
   
$\{\text{2,4,1}\}$ & $21.9603$   &  $\{\text{100,1,2}\}$ & $1.0495$     \\

$\{\text{}26,1,1\}$ & $18.7547$   &  $\{\text{55,2,1}\}$ & $0.9638$     \\

$\{\text{100,1,2}\}$ & $18.3988$   &  $\{\text{24,1,1}\}$ & $0.9383$     \\
    
$\{\text{115,1,2}\}$ & $16.9789$   &  $\{\text{26,1,1}\}$ & $0.8261$     \\

$\{\text{27,1,1}\}$ & $16.3683$  &  $\{\text{98,1,2}\}$ & $0.8062$     \\

$\{\text{6,1,1}\}$ & $15.1006$   &  $\{\text{131,2,2}\}$ & $0.7439$     \\

$\{\text{99,1,2}\}$ & $13.8501$   &  $\{\text{48,2,1}\}$ & $0.7243$      \\

$\{\text{172,3,2}\}$ & $12.7657$   &  $\{\text{176,3,2}\}$ & $0.6906$     \\

$\{\text{176,3,2}\}$ & $12.7579$   &  $\{\text{162,3,2}\}$ & $0.6773$     \\

$\{\text{13,1,1}\}$ & $12.6350$  &  $\{\text{155,3,2}\}$ & $0.6752$     \\

     \hline
\end{tabular}
}
\end{table}

\begin{table}
\caption{Multilayer subgraph centrality obtained with the modified tensor exponential 
(MSC$_{\exp_{0}}$) and multilayer subgraph centrality with the modified tensor resolvent 
(MSC$_{\rm res_0}$) for some nodes of the network in Example 6, with 
$\alpha=0.3/\lambda_{\max}$, $\beta=0.3$, $P=11$, and $m=5$.}

\label{tab:table11}
\centering
{\small
\begin{tabular}{ccc}
$\{i,\ell_1,\ell_2\}$ &  $MSC_{\exp_{0}}\{i,\ell_1,\ell_2\}$  & $MSC_{\rm res_0}\{i,\ell_1,\ell_2\}$  \\
   \hline
   
$\{\text{2,4,1}\}$ & $2.5074$     & $1.0082$     \\

$\{\text{}26,1,1\}$ & $2.5074$    & $1.0082$     \\

$\{\text{100,1,2}\}$ & $1.0000$  & $1.0000$     \\
    
$\{\text{115,1,2}\}$ & $1.0000$   & $1.0000$     \\

$\{\text{27,1,1}\}$ & $2.5074$   & $1.0082$     \\

$\{\text{6,1,1}\}$ & $2.5074$    & $1.0082$     \\

$\{\text{99,1,2}\}$ & $1.0000$   & $1.0000$      \\

$\{\text{172,3,2}\}$ & $1.0000$    & $1.0000$     \\

$\{\text{176,3,2}\}$ & $1.0914$    & $1.0027$     \\

$\{\text{13,1,1}\}$ & $1.0000$   & $1.0000$     \\

     \hline
\end{tabular}
}
\end{table}

\section{Conclusion}\label{sec6}
This paper investigates centrality measures for multilayer networks by introducing the 
exponential and the resolvent of the adjacency tensor associated with this network 
using the Einstein product. We showed how to approximate these tensor functions 
via Krylov subspace methods based on the tensor format. Numerical tests gave
satisfactory results. The paper illustrates the tensors are useful for modeling
multilayer networks and can be used to evaluate small to quite large networks. However,
the computations for very large networks and may require the use of 
parallel computers with many processors. This will be explored in future work.

\section*{Acknowledgment}
The authors would like to thank the referees for comments that improved the presentation.
Research by SN was partially supported by a grant from SAPIENZA Universit\`a di Roma 
and by INdAM-GNCS.

\end{document}